%% file: agt-2-8.tex
\documentclass{gtart}

\input agtout

\lognumber{8}
\volumenumber{2}
\volumeyear{2002}
\papernumber{8}
\published{6 March 2002}
\pagenumbers{157}{170}
\received{12 December 2001}
\revised{24 February 2002}
\accepted{28 February 2002}

\usepackage{amsmath}
\usepackage{amssymb}
\usepackage{epsf}

\newcommand{\scut}{\mathbf{S}_{\fr{E}}}
\newcommand{\csys}{\lbrace \fr{a}_1, \dots ,\fr{a}_p,\fr{b}_{11}, \dots ,\fr{b}_{1q_1}, \dots , \fr{b}_{r1}, \dots ,\fr{b}_{rq_r}, \fr{c}_1, \dots , \fr{c}_s \rbrace}
\newcommand{\mtwist}{D_{\fr{a}_1}^{\alpha _1} \cdots D_{\fr{a}_p}^{\alpha_p}D_{\fr{b}_{11}}^{\beta_{11}} \cdots       D_{\fr{b}_{1q_1}}^{\beta _{1q_1}} \cdots D_{\fr{b}_{r1}}^{\beta_{r1}} \cdots D_{\fr{b}_{rq_r}}^{\beta_{rq_r}}D_{\fr{c}_1}^{\gamma_1} \cdots D_{\fr{c}_s}^{\gamma_s}}
\newcommand{\lai}{\langle}
\newcommand{\rai}{\rangle}
\newcommand{\fr}[1]{\mathfrak{#1}}

\newtheorem{lemma}{Lemma}[section]
\newtheorem{theorem}{Theorem}[section]
\newtheorem{corollary}{Corollary}[section]

\begin{document}
\title{Abelian Subgroups of the Torelli Group}
\author{William R. Vautaw}
\address{Department of Mathematics, Michigan State University\\East 
Lansing, MI 48824, USA}

\email{vautawwi@pilot.msu.edu}

\begin{abstract}
Let $\mathbf{S}$ be a closed oriented surface of genus $g \geq 2$, and let $\mathcal{T}$ denote its Torelli group.  First, given a set $\fr{E}$ of homotopically nontrivial, pairwise disjoint, pairwise nonisotopic simple closed curves on $\mathbf{S}$, we determine precisely when a multitwist on $\fr{E}$ is an element of $\mathcal{T}$ by defining an equivalence relation on $\fr{E}$ and then applying graph theory.  Second, we prove that an arbitrary Abelian subgroup of $\mathcal{T}$ has rank $\leq 2g-3$.
\end{abstract}
\asciiabstract{
Let S be a closed oriented surface of genus g > 1, and let T denote
its Torelli group.  First, given a set E of homotopically nontrivial,
pairwise disjoint, pairwise nonisotopic simple closed curves on S, we
determine precisely when a multitwist on E is an element of T by
defining an equivalence relation on E and then applying graph theory.
Second, we prove that an arbitrary Abelian subgroup of T has rank <
2g-4.}

\maketitle

\section{Introduction}

Here we present the notation, definitions, and terminology that will be used in the paper.

\subsection{Surfaces}

Throughout this work, $\mathbf{S}$ will denote a closed, connected, oriented surface.  We use the symbols $\fr{a,b,c,e,h}$ to denote simple closed curves on $\mathbf{S}$. 

The \emph{mapping class group}, $\mathcal{M}(\mathbf{S})$, of $\mathbf{S}$ is the group of isotopy classes of  orientation preserving self-homeomorphisms of $\mathbf{S}$.  In general, we will not distinguish between a map $f \co \mathbf{S} \to \mathbf{S}$ and its isotopy class.  The symbol $D_\fr{e}$ will denote the right Dehn twist about the simple closed curve $\fr{e}$.  Recall that if $\fr{a}$ and $\fr{b}$ are simple closed oriented curves on $\mathbf{S}$, then in $H_1(\mathbf{S})$, the first homology group of $\mathbf{S}$ with integer coefficients, we have \[ D_\fr{a}(\fr{b}) = \fr{b} + \lai \fr{a},\fr{b} \rai \fr{a} \] where $\lai \fr{a},\fr{b} \rai$ denotes the algebraic intersection number of $\fr{a}$ and $\fr{b}$.  Also, the Dehn twists $D_{\fr{a}_1}$ and $D_{\fr{a}_2}$ commute if and only if the isotopy classes of the curves $\fr{a}_1$ and $\fr{a}_2$ have representatives that are disjoint.

The \emph{Torelli group}, $\mathcal{T}=\mathcal{T}(\mathbf{S})$, of $\mathbf{S}$ is the subgroup of the mapping class group consisting of the isotopy classes of those self-homeomorphisms of $\mathbf{S}$ which induce the identity isomorphism  on $H_1(\mathbf{S})$.  The Torelli group is torsion-free, and is trivial in the case of the sphere or torus.

\subsection{Graphs} 

We use graph-theoretic terminology consistent with its use in \cite{2}.  We remind the reader of the less familiar terms, and give the graph-theoretic definitions of those terms that may be used in different ways in ordinary topology. 

Throughout this work, $G$ will denote a connected, finite linear graph.  We include the possibility that $G$ may contain loops or parallel edges.  $E=E(G)$ will denote the edge set of $G$, and we use the symbols $a,b,c,e$ to denote edges of $G$.  For $E' \subset E(G)$, $G-E'$ denotes the subgraph obtained from $G$ by deleting the edges in $E'$, while $G+E''$ is the graph obtained from $G$ by adding a set of edges $E''$.    If $E= \lbrace e \rbrace$, then we write $G-e$ and $G+e$ instead of $G-\lbrace e \rbrace$ and $G+\lbrace e \rbrace$.  A \emph{bond} $E'$ in $G$ is a minimal subset of $E(G)$ such that $G-E'$ is disconnected. Note that $G-E'$ consists of precisely two components.  We say that the edge $e$ is a \emph{cut edge} if $G-e$ is disconnected.  We use the symbols $u,v,x,y$ to denote vertices of $G$.  The \emph{degree} of a vertex $v$ is the number of edges incident with $v$, each loop counting as two edges.

A $(v_0,v_n)$\emph{--walk} $W$ of \emph{length} $n$ is a finite nonempty alternating sequence, $W=v_0e_1v_1e_2v_2 \dots e_nv_n$, of vertices and edges such that the ends of the edge $e_i$ are the vertices $v_{i-1}$ and $v_{i}$ for $1 \leq i \leq n$.  If the edges of $W$ are distinct, $W$ is called a \emph{trail}.  A \emph{cycle} in $G$ is a closed trail of positive length whose origin and internal vertices are distinct.  Thus a cycle is an embedded circle in $G$.  For our purposes, to denote a trail or cycle, it will be enough to give its sequence of edges, and we do not distinguish between a closed trail $W$ and another closed trail whose sequence of edges is a cyclic permutation of $W$'s.
 
A \emph{spanning tree} $T$ is a subgraph of $G$ with the same vertex set as $G$ such that $T$ contains no cycles.  The number of edges in any spanning tree is equal to one less than the number of vertices of $G$.  Note that if $T$ is a spanning tree, and $e$ is an edge of $G$ not in $T$, then $T+e$ contains a unique cycle $C$, and $e$ is an edge of $C$, so the rank of $\pi_1(G)$ is equal to the number of edges of $G$ outside any spanning tree.  Every connected graph contains a spanning tree.

\section{Reduction Systems and Reduction System Graphs}

By a \emph{reduction system} $\fr{E}$ on $\mathbf{S}$ we mean a collection of simple closed curves on $\mathbf{S}$ that are homotopically nontrivial, pairwise disjoint, and pairwise nonisotopic. We use the symbols $\fr{a,b,c,e}$ to denote the elements of a reduction system $\fr{E}$, and $\mathbf{S}_\fr{E}$ to denote the natural compactification of $\mathbf{S}\diagdown \fr{E}$; that is,  ``$\mathbf{S}$ cut along $\fr{E}$.''  

We partition the set $\fr{E}=\lbrace \fr{e}_1,\fr{e}_2, \dots , \fr{e}_n \rbrace$ according to the equivalence relation $\sim$ generated by the rule
\begin{displaymath}
\fr{e}_i \sim \fr{e}_j\quad  \mathrm{if}\quad  \left\{ \begin{array}{ll}
\fr{e}_i = \fr{e}_j\\
\mathrm{or}\\
\lbrace \fr{e}_i,\fr{e}_j \rbrace \mathrm{\ is\ a\ minimal\ separating\ set\ in\ }\fr{E}.
\end{array} \right.
\end{displaymath}
Here, ``$\lbrace \fr{e}_i,\fr{e}_j \rbrace$ is a minimal separating set'' means that $\mathbf{S}_{\{ \fr{e}_i,\fr{e}_j \}}$ is disconnected, but both $\mathbf{S}_{\{\fr{e}_i \}}$ and $\mathbf{S}_{\{\fr{e}_j\}}$ are connected.  There are three types of $\sim$--equivalence classes:
\begin{itemize}
\item[\rm(i)]  Singleton classes $\{\fr{a}_1\}, \{\fr{a}_2\}, \dots , \{\fr{a}_p\}$ consisting of the separating curves $\fr{a}_1, \fr{a}_2, \dots , \fr{a}_p$ in $\fr{E}$. Such a curve wil be called an \emph{a--type curve}. 
\item[\rm(ii)] Classes $\{\fr{b}_{11}, \dots ,\fr{b}_{1q_1}\}, \{\fr{b}_{21}, \dots ,\fr{b}_{2q_2}\}, \dots , \{\fr{b}_{r1}, \dots ,\fr{b}_{rq_r}\}$ of cardinality\break at least 2.  Each such class $\{\fr{b}_{i1}, \dots , \fr{b}_{in_i}\}$ is characterized by the following three properties:
\begin{itemize}
\item[\rm(a)] No curve $\fr{b}_{ij}$ is separating.
\item[\rm(b)] $\fr{b}_{ij}$ is homologous to $\fr{b}_{ij'}$ for every pair $\fr{b}_{ij}$, $\fr{b}_{ij'}$.
\item[\rm(c)] Maximal with respect to (a) and (b).
\end{itemize}
A curve in such a class will be called a \emph{b--type curve}.
\item[\rm(iii)] Singleton classes $\{\fr{c}_1\},\{\fr{c}_2\}, \dots , \{\fr{c}_s\}$ where each $\fr{c}_i$ is non-separating and is homologous to no other curve in $\fr{E}$.  Such a curve will be called a \emph{c--type curve}.
\end{itemize}
According to (i), (ii), and (iii) above, we write \[ \fr{E}=\csys. \]
We use $\fr{E}$ to define a graph $G_\fr{E}$, which we call the \emph{reduction system graph} of $\fr{E}$, as follows:
\begin{itemize}
\item The vertices of $G_\fr{E}$ correspond to the components of $\scut$.
\item The edges of $G_\fr{E}$ correspond to the curves in the reduction system $\fr{E}$, with:
\begin{itemize}
\item (Links)  Two distinct vertices are connected by the edge $e_i$ if and only if the curve $\fr{e}_i$ in $\fr{E}$ is a common boundary curve of the two components of $\scut$ which correspond to the vertices in question.
\item (Loops)  A vertex has a loop $e_i$ if and only if the curve $\fr{e}_i$ in $\fr{E}$ represents two boundary curves of the component of $\scut$ which corresponds to the vertex in question.
\end{itemize}
\end{itemize}
Note that $G_\fr{E}$ is connected, and that any connected graph $G$ is $G_\fr{E}$ for some surface $\mathbf{S}$ and some reduction system $\fr{E}$ on $\mathbf{S}$.  However, the genus of $\mathbf{S}$ is not determined by $G$, any two possible $\mathbf{S}$'s differing by the genera of their complementary components.  But, unless $G$ is the graph consisting of a single vertex and either no edges or a single loop, then $\mathrm{genus}(\mathbf{S}) \geq \mathrm{rank}(\pi_1(G))+(\mathrm{number\ of\ vertices\ of\ degree }\leq 2).$

Since $\mathbf{S}$ and $\fr{E}$ will be fixed, we will denote $G_\fr{E}$ simply by $G$.

The $\sim$--equivalence relation on the curves in $\fr{E}$ induces a $\sim$--equivalence relation on the edge set $E(G)=\{e_1,e_2,...,e_n\}$ of $G$.  It is generated by 
\begin{displaymath}
e_i \sim e_j\quad  \mathrm{if}\quad  \left\{ \begin{array}{ll}
e_i = e_j\\
\mathrm{or}\\
\lbrace e_i,e_j \rbrace \mathrm{\ is\ a\ bond.}
\end{array} \right.
\end{displaymath}
(Again, it should be noted that this equivalence relation may be defined for $any$ graph $G$.)  The three types of equivalence classes described above become, for $G$, 
\begin{itemize}
\item[\rm(i)] Singleton classes $\{a_1\},\dots ,\{a_p\}$ consisting of the cut edges $a_1, \dots ,a_p$ of $G$. Such an edge will be called an \emph{a-type edge}.
\item[\rm(ii)] Classes $\{b_{11}, \dots , b_{1q_1}\}, \{b_{21}, \dots , b_{2q_2}\}, \dots ,\{b_{r1}, \dots ,b_{rq_r}\}$ of cardinality at least 2.  Each such class is characterized by the following three properties:
\begin{itemize}
\item[\rm(a)] No edge $b_{ij}$ is a cut edge. 
\item[\rm(b)] $\{b_{ij},b_{ij'} \}$ is a bond for every pair $b_{ij}$, $b_{ij'}$.
\item[\rm(c)]  Maximal with respect to (a) and (b). 
\end{itemize}
An edge in such a class will be called a \emph{b-type edge}.
\item[\rm(iii)]  Singleton classes $\{c_1\}, \dots ,\{c_s\}$ where each $c_i$ is not a cut edge, and forms a 2--edge bond with no other edge of $G$.  Such an edge will be called a \emph{c-type edge}.
\end{itemize}
According to (i), (ii), and (iii) above, we write \[E(G)~=~\lbrace a_1, \dots , a_p, b_{11}, \dots , b_{1q_1}, b_{21},\dots ,b_{2q_2}, \dots , b_{r1}, \dots , b_{rq_r}, c_1, \dots ,c_s \rbrace.\]

\begin{figure}[ht!]
\cl{\epsfxsize=2.5in  \epsfbox{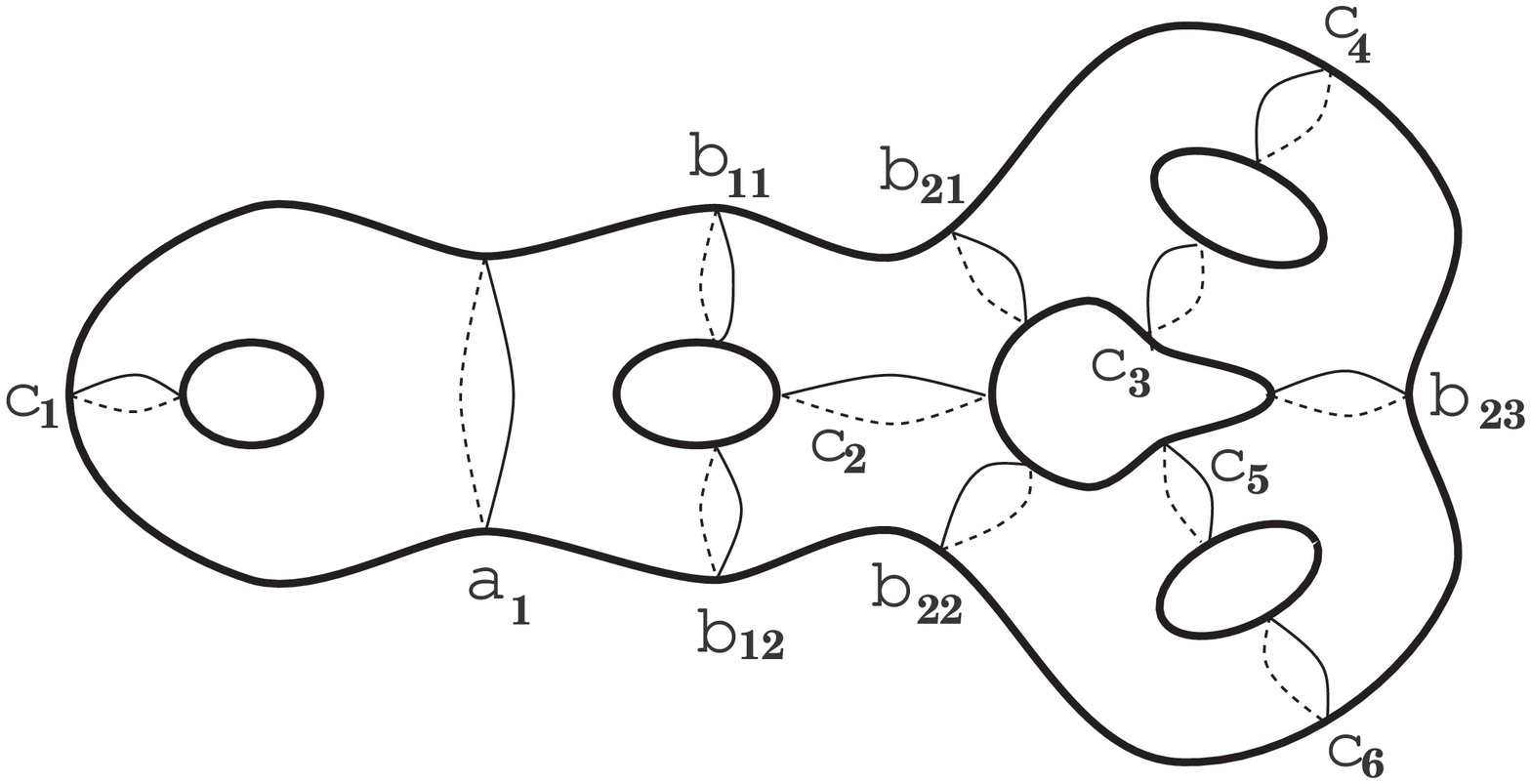} \epsfxsize=2.5in \epsfbox{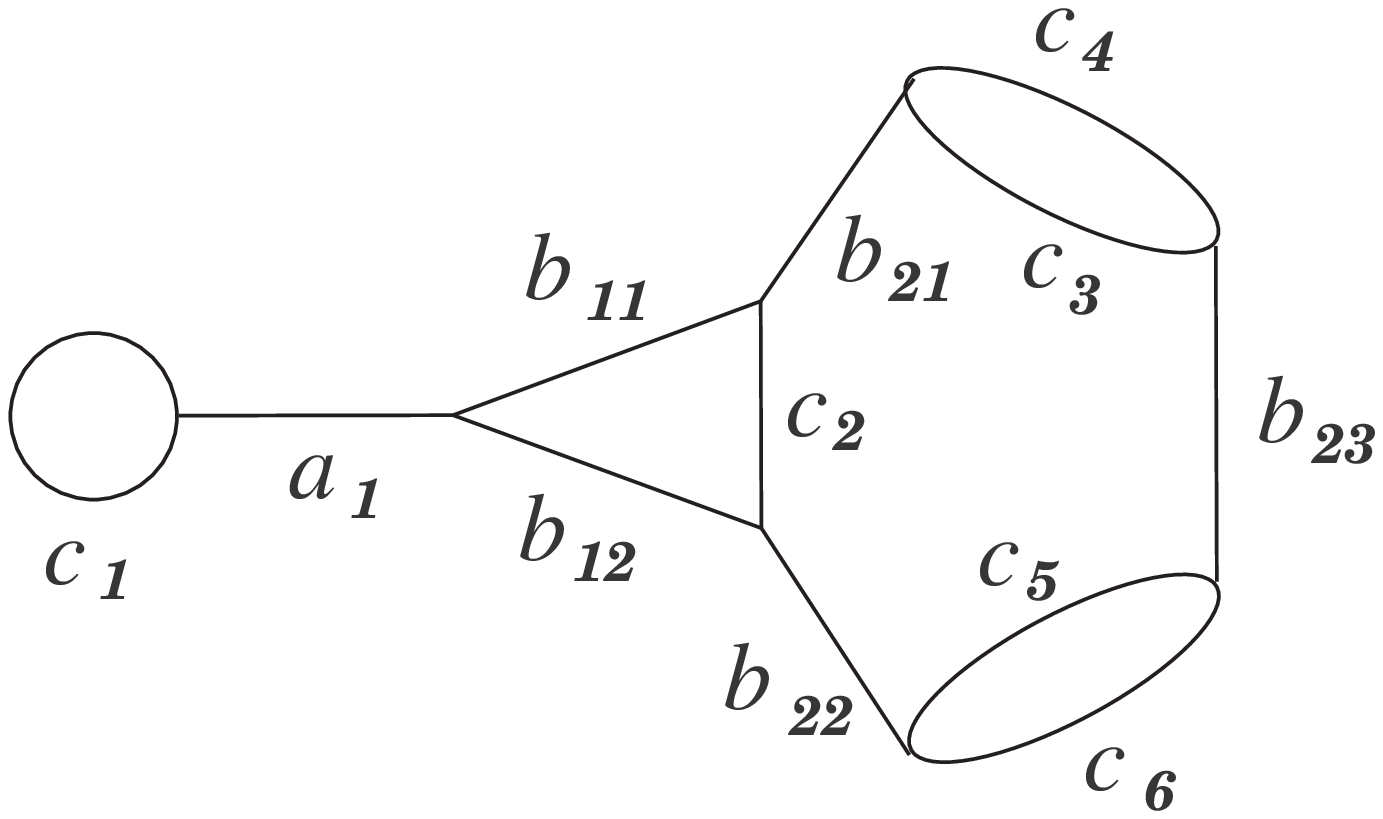}}
\vglue 3pt
\end{figure}

A typical example is shown above.

Now let $\fr{h}$ be a simple closed curve on $\mathbf{S}$ that intersects each element of $\fr{E}$ transversely at most once.  Starting at any point on $\fr{h}$ and travelling in either direction gives a cyclic ordering of the reduction curves which $\fr{h}$ intersects, thus defining a closed trail $H$ in $G$.    Note that $H$ is a cycle in $G$ if and only if $\fr{h} \cap \mathbf{S}_i$ is either empty or is a single (that is, connected) arc, for every component $\mathbf{S}_i$ of $\scut$.  Likewise, given a closed trail $H$ in $G$, there is such a curve $\fr{h}$ on $\mathbf{S}$ defining $H$.  The fact that the isotopy class of $\fr{h}$ is never unique is not important for our purposes.
  
The following figure shows a typical example.  Note that $\fr{h}_1$ and $\fr{h}_2$ are nonisotopic curves which both define the cycle $H = b_{11}c_2b_{12}$.

\begin{figure}[ht!]
\cl{\epsfxsize=2.5in  \epsfbox{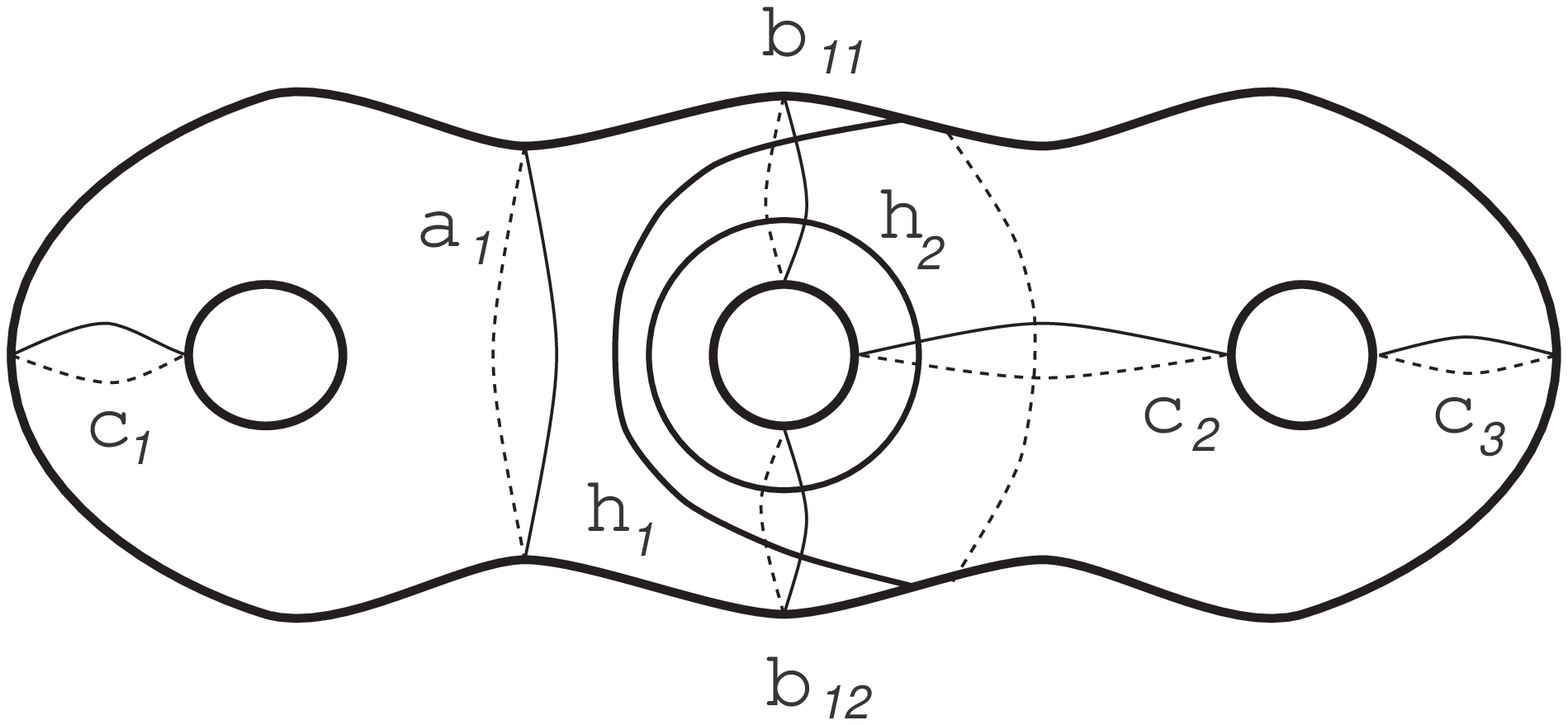} \epsfxsize=2.5in
\epsfbox{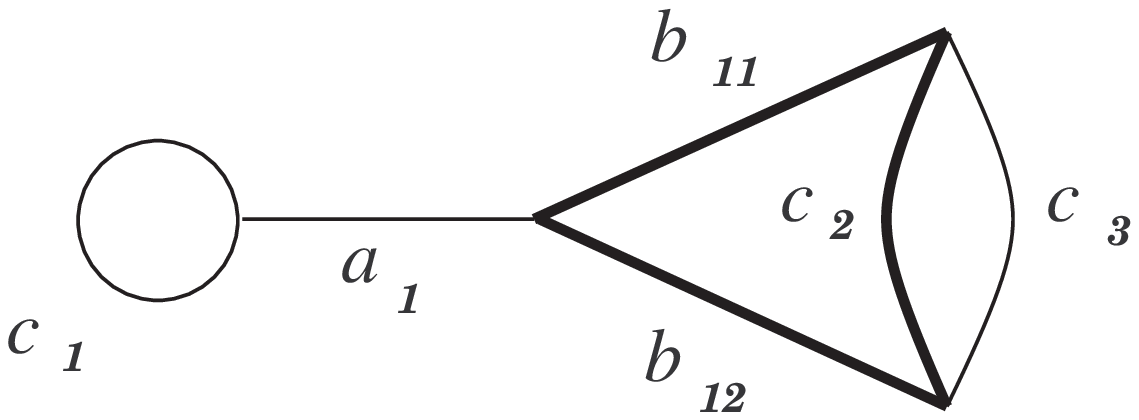}}
\vglue 3pt
\end{figure}

The remainder of this section presents some purely graph-theoretic results, concluding with Theorem 2.1, which is used in the following section.  So for the remainder of this section, let $G$ denote an arbitrary connected graph.  We explain here the notation and terminology we use.  Given a subgraph $H$ of $G$, we let $G\bullet H$ denote the graph obtained by deleting every edge $e$ of $H$ and identifying the ends of $e$.  Equivalently, thinking of $G$ as a CW--complex and $H$ as a subcomplex, $G\bullet H$ is the complex obtained from $G$ by crushing each component of $H$ to a point.  Thus, we have a quotient (``contraction'') map $p \co G \to G\bullet H$.  Next, by a \emph{cut vertex} of $G$, we mean a vertex $v$ of $G$ such that when $v$, \emph{and only} $v$, is removed from the topological space $G$, the resulting space is disconnected.  (This is not the definition used by graph theorists, but is an equivalent topological one.)  A \emph{block} is a connected graph without cut vertices, and a \emph{block of a graph} is a subgraph that is a block and is maximal with repsect to that property.  Any graph is the union of its blocks.  We leave the proofs of the first two lemmas to the reader.

\begin{lemma}
If $G$ has no cut edges, then any two vertices of $G$ are connected by two edge-disjoint paths.
\end{lemma}

\begin{lemma}
Let $b_1$ and $b_2$ be edges of $G$ such that $\lbrace b_1, b_2 \rbrace$ is a bond.  If $C$ is a cycle in $G$, and $b_1$ is an edge of $C$, then so is $b_2$.
\end{lemma}

\begin{lemma}
Let $c$ be a \textrm{c}--type edge in $G$ that is not a loop.  Then $c$ is contained in two cycles, the intersection of whose edge sets is precisely $c$.
\end{lemma}

\begin{proof}
Assume that $G$ is a block.  If $G$ has exactly two vertices, then each edge of $G$ is a link, and $G$ must have at least three edges, since $c$ is a c--type edge.  The result is clear in this case.  Otherwise, $G$ has at least three vertices and no cut edges.  Consider the graph $G-c$.  If $G-c$ has a cut edge $e$, then $G - \lbrace c,e \rbrace$ is not connected, so $\lbrace c,e \rbrace$ is a bond of $G$.  This contradicts the fact that $c$ is a c--type edge.  So $G-c$ has no cut edges.  By Lemma~2.1, there are two edge-disjoint paths $P$ and $P'$ in $G-c$ connecting the ends of $c$.  Then the cycles $C=P+c$ and $C'=P'+c$ have exactly the edge $c$ in common.  In the case that $G$ is not a block, we let $B$ be the block of $G$ containing $c$.  It is easy to see that $c$ is a c--type edge of $B$, so we apply the first case to $B$ and find two such cycles within $B$. 
\end{proof}

\begin{theorem}
Let $G$ have edge set
$$E(G)=\{ a_1, \dots , a_p, b_{11}, \dots , b_{1q_1}, \dots , b_{r1}, \dots , b_{rq_r}, c_1, \dots ,c_s \},$$ notated according to $\mathrm{a}$--~, $\mathrm{b}$--, and $\mathrm{c}$--type equivalence classes.  Let $w \co E(G)\to \mathbb{Z}$ be a weighting of $G$.  Then $w(H)=0$ for every cycle $H$ in $G$ if and only if
\begin{itemize}
\item[\rm(i)] $w(c_i) = 0$,\ \  $1 \leq i \leq s$,\ \ and
\item[\rm(ii)] $w(b_{j1}) + w(b_{j2}) + \cdots + w(b_{jq_j}) = 0$,\ \  $1 \leq j \leq r$.
\end{itemize}
\end{theorem}

\begin{proof}
$\Rightarrow$\qua Assume that $w(H)=0$ for every cycle $H$ in $G$.

(i)\qua  Let $c$ be a c--type edge with ends $u$ and $v$.  If $c$ is a loop, then $w(c)=0$, by hypothesis.  Otherwise, there are two edge-disjoint $(u,v)$--paths, $P$ and $\overline{P'}$, in $G-c$.  We have three cycles: $P+c$, $P'+c$, and $P+P'$.  Thus:
\begin{displaymath}
\left.
\begin{array}{rr}
w(P)+w(c)=w(P+c) = 0\\ 
w(P')+w(c)=w(P'+c) = 0\\
w(P)+w(P')=w(P+P') = 0
\end{array} \right\}
\Longrightarrow w(c) = 0
\end{displaymath}
(ii)\qua  Let $B$ be the equivalence class of the b--type edge $b$, and $\overline{B} = B \diagdown \lbrace b \rbrace$.  Let $p \co G \to G \bullet \overline{B}$ be the contraction map.  Suppose that $b$ is a cut edge of $G\bullet \overline{B}$, separating it into two components $G_1$ and $G_2$.  Then the restriction of $p$ to $G-b$ maps onto the disconnected space $G_1 \cup G_2$, and so $G-b$ is disconnected.  This is a contradiction to the hypothesis that $b$ is a b--type edge of $G$.  We obtain a similar contradiction if we suppose $\{b,e\}$ is a bond in $G\bullet \overline{B}$.  Thus $b$ is a c--type edge in $G \bullet \overline{B}$.  If $b$ is a loop in $G \bullet \overline{B}$, then $p^{-1}(b)= B$, which therefore forms a cycle in $G$.  So equation (ii) holds for the equivalence class of $b$.

If $b$ is not a loop in $G \bullet \overline{B}$, then by Lemma 2.3 there are two cycles $\overline{H}$ and $\overline{H'}$ in $G \bullet \overline
{B}$, the intersection of whose edge sets is $\{b\}$.  Lemma 2.2 implies that $p^{-1}(\overline{H})$ and $p^{-1}(\overline{H'})$ are cycles $H$ and $H'$, respectively, the intersection of whose edge sets is precisely $B$. Thus we have
\begin{displaymath}
\left. 
\begin{array}{rr}
0=w(H)=w(B)+w(H-B)\\ 
0=w(H')=w(B)+w(H'-B)
\end{array} \right\}
\Longrightarrow\ 0=2w(B)+w(H\Delta H')=2w(B).
\end{displaymath}
And so, $w(B)=0$.  Here we have used the fact that the symmetric difference $H\Delta H'$ of the cycles $H$ and $H'$ is a disjoint union of cycles (regarded as sets of edges).

$\Leftarrow$\qua  Assume that
\begin{itemize}
\item[\rm(i)] $w(c_i) = 0$,\ \  $1 \leq i \leq s$,\ \ and
\item[\rm(ii)] $w(b_{j1}) + w(b_{j2}) + \cdots + w(b_{jq_j}) = 0$,\ \  $1 \leq j \leq r$.
\end{itemize}
Let $H$ be a cycle in $G$.  $H$ contains no a--type edges, since they are cut edges, and by Lemma~2.2, if $H$ contains one edge of a b--type class, then it contains the whole class.  So the assumptions imply that $w(H)=0$.
\end{proof}

\section{Abelian Subgroups in the Torelli Group}

We at first consider a specific type of Abelian subgroup of the Torelli group $\mathcal{T}(\mathbf{S})$, namely one consisting of multitwists --- that is, compositions of left and right Dehn twists about a fixed reduction system on $\mathbf{S}$.

\begin{theorem}
Let $\mathbf{S}$ be a closed, connected, oriented surface, and let \[ \fr{E}=\csys \] be a reduction system on $\mathbf{S}$, notated by $\mathrm{a}$--, $\mathrm{b}$--, and $\mathrm{c}$--type $\sim$--equivalence classes as in section 2.  Let $\mathcal{D}_\fr{E}$ be the multitwist group on $\fr{E}$, and let \[ f=\mtwist \] be an element of $\mathcal{D}_\fr{E}$.  Then $f$ is an element of $\mathcal{D}_\fr{E} \cap \mathcal{T} \equiv \mathcal{T}_{\fr{E}}$, which we call the \emph {Torelli multitwist group of $\fr{E}$}, if and only if
\begin{itemize}
\item[\rm(i)] $\gamma _i = 0$,\  $1 \leq i \leq s$, and
\item[\rm(ii)] $\beta_{j1}+\beta_{j2}+\cdots +\beta_{jq_j} = 0$,\  $1 \leq j \leq r$.
\end{itemize}
\noindent Consequently, $\mathcal{T}_{\fr{E}}$ is a free Abelian group of rank \[ p+(q_1-1)+(q_2-1)+ \cdots +(q_r-1) = p+q_1+q_2+ \cdots +q_r-r. \]
\end{theorem}

\rk{Remark}  A set of $\sim$--equivalence class representatives of the curves in $\fr{E}$ is in general \emph{not} linearly independent in $H_1(\mathbf{S})$, so the nondegeneracy of the algebraic intersection $\lai\ ,\ \rai$ is not sufficient to prove the theorem.

\proof
\noindent $\Rightarrow$\qua  Assume that $f \in \mathcal{T}_\fr{E}$.

Let $G$ be the reduction system graph of $\fr{E}$ with edge set $E(G)$.  We weight each edge of $G$ according to the exponent in $f$ of the twist about its corresponding curve in $\fr{E}$, giving $w \co E(G)\to \mathbb{Z}$. 

Let $H=e_1e_2,\dots ,e_n$ be a cycle in $G$.  Then, as in section 2, $H$ is defined by any simple closed curve $\fr{h}$ on $\mathbf{S}$ that intersects each of the corresponding curves $\fr{e}_1,\fr{e}_2, \dots,\fr{e}_n$ of $\fr{E}$ exactly once, and does not intersect any of the other curves of $\fr{E}$.  Orient $\fr{h}$.  Then orient the curves $\fr{e}_1,\fr{e}_2, \dots,\fr{e}_n$ so that $\lai \fr{e}_i,\fr{h} \rai=1$. So we have \[0 =\lai \fr{h},\fr{h} \rai = \lai \fr{h},f(\fr{h})\rai = \lai \fr{h},\fr{h}+ \epsilon_1\fr{e}_1+\epsilon_2\fr{e}_2+ \cdots +\epsilon_n\fr{e}_n\rai=\epsilon_1 +\epsilon_2+\cdots +\epsilon_n,\]  
where $\epsilon_i=w(e_i)$.  Hence the weight of every cycle in $G$ is zero.  The conclusion follows from Theorem 2.1.  

\noindent $\Leftarrow$\qua  Assume that 
\begin{itemize}
\item[\rm(i)] $\gamma _i = 0$,\ $1 \leq i \leq s$,\ \  and 
\item[\rm(ii)] $\beta_{j1}+\beta_{j2}+\cdots +\beta_{jq_j} = 0$,\ $1 \leq j \leq r$.
\end{itemize}

Since $H_1(\mathbf{S})$ has a basis consisting of simple closed curves, in order to prove that $f \in \mathcal{T}$, it suffices to show that in $H_1(\mathbf{S})$, we have $f(\fr{h})=\fr{h}$ for any simple closed curve $\fr{h}$ on $\mathbf{S}$.  Note that for any such $\fr{h}$, we have $\lai \fr{a}_i, \fr{h} \rai =0$, $1 \leq i \leq p$, and after orienting $\fr{h}$ and then each $\fr{b}_{ij}$ so that $\lai \fr{b}_{ij},\fr{h} \rai = \lai \fr{b}_{i1}, \fr{h} \rai$, we have $\fr{b}_{ij}=\fr{b}_{i1}$, $2 \leq j \leq q_i$, $1 \leq i \leq r$.  Let $\delta_i =\lai \fr{b}_{i1},\fr{h} \rai$. Then in $H_1(\mathbf{S})$ we have: 
\begin{align}
f(\fr{h}) & = \mtwist (\fr{h})\notag\\
 & = \fr{h}+\beta_{11}\lai \fr{b}_{11},\fr{h}\rai \fr{b}_{11}+\cdots+\beta_{1q_1}\lai \fr{b}_{1q_1},\fr{h}\rai \fr{b}_{1q_1}+\cdots\notag\\
  & \ \ \,\quad  +\beta_{r1}\lai \fr{b}_{r1},\fr{h}\rai \fr{b}_{r1}+\cdots+\beta_{rq_r}\lai \fr{b}_{rq_r},\fr{h}\rai \fr{b}_{rq_r}\notag\\
& = \fr{h} +\delta_1(\beta_{11}+ \cdots +\beta_{1q_1})\fr{b}_{11}+ \cdots +\delta_r(\beta_{r1}+ \cdots +\beta_{rq_r})\fr{b}_{r1}\notag\\
& = \fr{h}\tag*{\qed}
\end{align}

\begin{theorem}
Let $\mathbf{S}$ be a closed connected oriented surface, and let $\fr{E}=\{\fr{e}_1,\fr{e}_2, \dots,\fr{e}_n\}$ be a reduction system on $\mathbf{S}$.  Let $f=D_{\fr{e}_1}^{\epsilon_1}D_{\fr{e}_2}^{\epsilon_2}\cdots D_{\fr{e}_n}^{\epsilon_n}$ be a multitwist on $\fr{E}$. Let $G$ be the reduction system graph of $\fr{E}$, and define a weighting $w \co E(G)\to \mathbb{Z}$ of $G$ by $w(e_i) = \epsilon_i$.  Then $f$ is in the Torelli multitwist group $\mathcal{T}_\fr{E}$ if and only if the weight of every cycle in $G$ is zero.
\end{theorem}
\begin{proof}
Partition $\fr{E}$ into $\sim$--equivalence classes, so \[\fr{E}= \csys .\]  Theorems 2.1 and 3.1 show the conditions to be equivalent.
\end{proof}

Given a pair, $\fr{e}_1$ and $\fr{e}_2$, of disjoint, non-separating, but homologous simple closed curves on $\mathbf{S}$, we call $D_{\fr{e}_1}D_{\fr{e}_2}^{-1}$ a \emph{bounding-pair map} or \emph{BP map}.  Powell \cite{5} has shown that the Torelli group $\mathcal{T}$ is generated by BP maps and Dehn twists about separating simple closed curves.

\begin{corollary}
Let $\mathbf{S}$, $\fr{E}$, $\mathcal{D}_\fr{E}$, and $\mathcal{T}_\fr{E}$ be as in Theorem~3.1.  Let $\mathcal{D'}$ be the subgroup of $\mathcal{M}(\mathbf{S})$ generated by 
\begin{itemize}
\item[\rm(i)] BP maps about bounding pairs in $\fr{E}$, and
\item[\rm(ii)] Dehn twists about separating curves in $\fr{E}$.
\end{itemize}
Then $\mathcal{D'} = \mathcal{D}_\fr{E} \cap \mathcal{T} = \mathcal{T}_\fr{E}$.
\end{corollary}

\begin{proof}
By the definition of $\mathcal{D}_\fr{E}$, it is clear that every generator of $\mathcal{D'}$ is in $\mathcal{D}_\fr{E}$.  By Powell's result noted above, every generator of $\mathcal{D'}$ is in $\mathcal{T}$.  Thus $\mathcal{D'} \subseteq \mathcal{D}_\fr{E} \cap \mathcal{T}$.  We must show that $\mathcal{D}_\fr{E} \cap \mathcal{T} \subseteq \mathcal{D'}$.

Let $f \in  \mathcal{D}_\fr{E} \cap \mathcal{T}$.  By Theorem~3.1, we know that 
\[f=D_{\fr{a}_1}^{\alpha _1} \cdots D_{\fr{a}_p}^{\alpha _p}D_{\fr{b}_{11}}^{\beta _{11}} \cdots D_{\fr{b}_{1q_1}}^{\beta _{1q_1}}D_{\fr{b}_{21}}^{\beta _{21}} \cdots D_{\fr{b}_{2q_2}}^{\beta _{2q_2}} \cdots D_{\fr{b}_{r1}}^{\beta_{r1}} \cdots D_{\fr{b}_{rq_r}}^{\beta _{rq_r}},\]
where $\beta_{i1}+\beta_{i2}+\cdots +\beta_{iq_i}=0$, \ $1 \leq i \leq r$.  Since each $D_{\fr{a}_i}^{\alpha_i}$ is a product of type--(ii) generators of $\mathcal{D'}$, we will be done if we write $D_{\fr{b}_{i1}}^{\beta _{i1}}D_{\fr{b}_{i2}}^{\beta _{i2}} \cdots D_{\fr{b}_{iq_i}}^{\beta _{iq_i}}$ as a product of BP maps.  We do this:  
\[ D_{\fr{b}_{i1}}^{\beta _{i1}}D_{\fr{b}_{i2}}^{\beta _{i2}} \cdots D_{\fr{b}_{iq_i}}^{\beta _{iq_i}}= (D_{\fr{b}_{i2}}D_{\fr{b}_{i1}}^{-1})^{\beta_{i2}}(D_{\fr{b}_{i3}}D_{\fr{b}_{i1}}^{-1})^{\beta_{i3}}\cdots (D_{\fr{b}_{iq_1}}D_{\fr{b}_{i1}}^{-1})^{\beta_{iq_i}}, \]
where we note that $-\beta_{i2}-\beta_{i3}-\cdots -\beta_{iq_i}=\beta_{i1}.$
\end{proof}

\begin{corollary}
Let $\mathbf{S}$ be a closed, connected, oriented surface, and let \[ \fr{E}=\csys \] be a reduction system on $\mathbf{S}$, notated by $\mathrm{a}$--, $\mathrm{b}$--, and $\mathrm{c}$--type $\sim$--equivalence classes as in section 2.  Let $\mathcal{D}_\fr{E}$ be the multitwist group on $\fr{E}$, and let \[ f=\mtwist \] be an element of $\mathcal{D}_\fr{E}$.  Let $m \geq 2$ be an integer.
 
Then $f \in \Gamma_{\mathbf{S}}(m) \equiv \{g \in \mathcal{M}(\mathbf{S}): g\ acts\ trivially\ on\ H_1(\mathbf{S};\mathbb{Z}_m)\}$ if and only if
\begin{itemize}
\item[\rm(i)] $\gamma _i \equiv 0\ (\mathrm{mod}\ m)$,\  $1 \leq i \leq s$, and
\item[\rm(ii)] $\beta_{j1}+\beta_{j2}+\cdots +\beta_{jq_j} \equiv 0\ (\mathrm{mod}\ m)$,\  $1 \leq j \leq r$.
\end{itemize}
\end{corollary}

Let $\mathbf{S}$ be the surface of genus $g \geq 2$ and $\fr{E}$ the reduction system on $\mathbf{S}$ shown below.  Since $\fr{E}$ consists of $2g-3$ $\mathrm{a}$--type curves, $\mathrm{rank}(\mathcal{T}_\fr{E}) =2g-3$.  This example, along with Theorem~4.1 below, shows that the maximal rank of an Abelian subgroup of the Torelli group is attained by a multitwist group.

\begin{figure}[ht!]
\cl{\epsfxsize=4in  \epsfbox{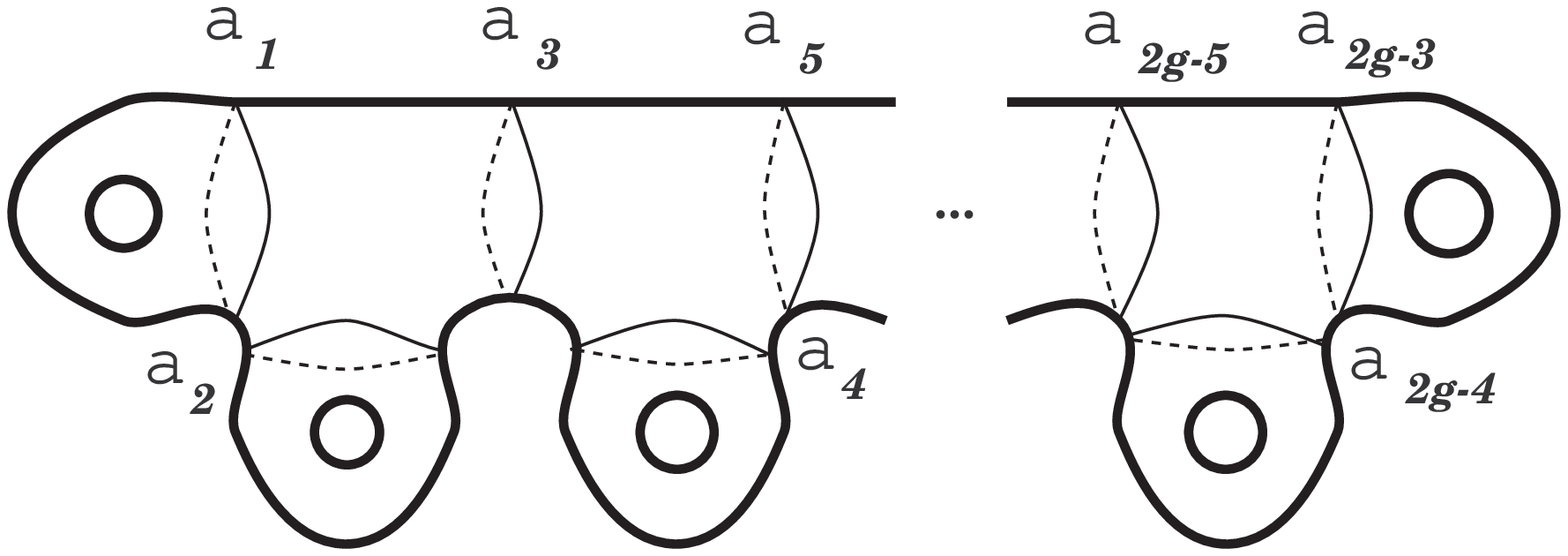}}
\vglue 3pt\end{figure}

\rk{Remark}  One particular naively-expected symplectic analogue of Theorem 3.1 is not true:

\textbf{``Conjecture''}\qua  Let $(\mathcal{V}, \lai,\rai)$ be a symplectic lattice of rank $2g$, where $g \geq 2$.  Let $\{v_1,v_2, \dots ,v_n \}$ be a set of primitive vectors in $\mathcal{V}$ that are pairwise linearly independent and symplectically orthogonal.  Let $T_i$ be the transvection corresponding to the vector $v_i$.  Thus $T_i(w)=w+\lai v_i, w \rai$ for any $w \in \mathcal{V}$.  Let $m_1,m_2, \dots ,m_n$ be integers.  Then the ``multitransvection'' $T\!=\!T_1^{m_1}T_2^{m_2} \cdots T_n^{m_n}$ is the identity on $\mathcal{V}$ if and only if $ m_i =0,\ 1 \leq i \leq n$.

But now let $\{a_1,b_1,a_2,b_2, \dots, a_g,b_g\}$ be the standard symplectic basis for $\mathcal{V}$, and for $i$ = 1, 2, 3, and 4, let $v_i=a_1+ia_2$.  Let $m_1=1$, $m_2=-3$, $m_3=3$, and $m_4=-1$.  One can verify that $T=T_1^{m_1}T_2^{m_2}T_3^{m_3}T_4^{m_4}=id_\mathcal{V}$.  This shows the conjecture to be false.

Now we prove that for any closed oriented surface of genus $g \geq 2$, the general Abelian subgoup of its Torelli group has rank $\leq 2g-3$.  We first give two lemmas.  

\begin{lemma}
Let $\mathbf{S}$ be a closed, connected, oriented surface, and $\fr{E}$ a reduction system on $\mathbf{S}$ with reduction system graph $G$.  Let $\mathcal{T}_\fr{E}$ be the Torelli multitwist group on $\fr{E}$, as in Theorem~3.1.  Then $\mathrm{rank}(\mathcal{T}_\fr{E}) \leq \nu-1$, where $\nu$ is the number of vertices of $G$, or, equivalently, the number of components of $\scut$.
\end{lemma}

\proof
Let $G$ have edge setl
$$E(G)=\lbrace a_1, \dots ,a_p,b_{11}, \dots ,b_{1q_1},b_{21}, \dots ,b_{2q_2}, \dots ,b_{r1}, \dots ,b_{rq_r},c_1, \dots , c_s \rbrace.$$  
Let $E'=\lbrace b_{11}, \dots ,b_{1(q_1-1)},b_{21},\dots ,b_{2(q_2-1)}, \dots ,b_{r1},\dots ,b_{r(q_r-1)} \rbrace \subseteq E(G)$, and let $G'=G[E']$.  Then $G'$ contains no cycles, since any cycle containing one edge of a b--type class contains the whole class.  Therefore, $G'$ is contained in a spanning tree $T$ of $G$.  Since each $a_i$ is a cut edge, $T$ contains $a_i$, $1 \leq i \leq p$.  

So $T$ contains the set of edges $E' \cup \lbrace a_1,a_2,\dots,a_p \rbrace$.  But by Theorem~3.1, the cardinality of this set is equal to the rank of $\mathcal{T}_\fr{E}$.  This gives us $$\nu -1 =\mathrm{card}(E(T)) \geq p+(q_1-1)+\cdots +(q_r-1) = \mathrm{rank}(\mathcal{T}_\fr{E}).\eqno{\qed}$$

\begin{lemma}
Let $\mathbf{S}$ be a closed, connected, oriented surface of genus $g \geq 2$, and let $\fr{E}$ be a reduction system on $\mathbf{S}$.  Let $\Omega$ denote the number of components of $\scut$ not homeomorphic to a pair of pants or a one-holed torus.  Let $\mathcal{T}_\fr{E}$ be the Torelli multitwist group on $\fr{E}$.  Then $\mathrm{rank}(\mathcal{T}_\fr{E}) +\Omega \leq 2g-3$. 
\end{lemma}

\proof
Let $G$ be the reduction system graph of $\fr{E}$.  We use the following notation:
\begin{itemize}
\item $\Gamma$ is the maximum genus of any component of $\scut$.
\item $\Delta$ is the maximum degree of any vertex of $G$, or, equivalently, the maximum number of boundary curves of any component of $\scut$.
\item $\nu_b$ is the number of vertices of $G$ of degree $b$, or, equivalently, the number of components of $\scut$ with $b$ boundary curves.
\item $\nu_b^\gamma$ ($\nu_b^{\geq \gamma}$) is the number of components of $\scut$ of genus $\gamma$ ($\geq \gamma$) having $b$ boundary curves, or, equivalently, the number of vertices of $G$ of degree $b$ corresponding to a component of $\scut$ of genus $\gamma$ ($\geq \gamma$).
\end{itemize}
So we have: \[\nu_b =\sum_{\gamma=0}^\Gamma \nu_b^\gamma \qquad \mathrm{and} \qquad \nu =\sum_{b=1}^\Delta \nu_b \] 
But the assumption that each element of $\fr{E}$ is homotopically nontrivial means $\nu_1^0=0$, and the assumption that the elements of $\fr{E}$ are pairwise nonisotopic means $\nu_2^0 = 0$. So, in fact, $\nu =\nu_1^{\geq 1}+\nu_2^{\geq 1}+\nu_3 +\nu_4 +\cdots +\nu_\Delta$. Now, $\nu_1^1$ is the number of one-holed tori, and $\nu_3^0$ is the number of pairs of pants, so by the definition of $\Omega$, we have $\Omega = \nu_1^{\geq 2}+\nu_2^{\geq 1}+\nu_3^{\geq 1} + \nu_4 + \cdots + \nu_\Delta$.  Hence
\begin{eqnarray*}
2g-2 & = & -\chi (\mathbf{S})\\
& = & \sum_{ \mathrm{components}\ \mathbf{V} \ \mathrm{of}\ \scut} -\chi (\mathbf{V})\\
& = & \sum_{\gamma=1}^\Gamma (2\gamma -1)\nu_1^\gamma + \sum_{\gamma=1}^\Gamma (2\gamma)\nu_2^\gamma + \sum_{b=3}^\Delta \sum_{\gamma=0}^\Gamma (2\gamma+b-2)\nu_b^\gamma
\end{eqnarray*}
By Lemma~3.1, $\mathrm{rank}(\mathcal{T}_\fr{E}) \leq \nu -1$, so we have
{\setlength\arraycolsep{1pt}
\begin{align*}
\mathrm{rank}(\mathcal{T}_\fr{E})&+ \Omega  \leq \nu + \Omega -1\\
& = (\nu_1^{\geq 1} +\nu_1^{\geq 2}+\cdots +\nu_\Delta)+(\nu_1^{\geq 2}+\nu_2^{\geq 1} +\nu_3^{\geq 1}+\nu_4 +\cdots +\nu_\Delta)-1\\
& = [(\nu_1^1+2\nu_1^{\geq 2})+2\nu_2^{\geq 1}+(\nu_3^0+2\nu_3^{\geq 1})+2\nu_4+2\nu_5+\cdots +2\nu_\Delta]-1\\
& \leq \Big[\sum_{\gamma=1}^\Gamma (2\gamma-1)\nu_1^\gamma + \sum_{\gamma=1}^\Gamma (2\gamma)\nu_2^\gamma + \sum_{b=3}^\Delta \sum_{\gamma=0}^\Gamma (2\gamma+b-2)\nu_b^\gamma\Big]-1\\
& = -\chi(\mathbf{S})-1\\
& = 2g-3\tag*{\qed}
\end{align*}}

\begin{theorem}
Let $\mathbf{S}$ be a closed, connected, oriented surface of genus $g \geq 2$, and let $\mathcal{A}$ be an Abelian subgroup of $\mathcal{T}$, the Torelli group of $\mathbf{S}$.  Then $\mathrm{rank}(\mathcal{A}) \leq 2g-3$.
\end{theorem}

\proof
This proof is an adaptation of a analogous proof in \cite{1}.  That paper also introduces the reduction homomorphism and essential reduction system which we refer to here.

Let $f \in \mathcal{A}$, $f \neq 0$.  By Thurston's classification, $f$ is either reducible, pseudo-Anosov, or of finite order.  Since $\mathcal{T}$ is torsion-free, $f$ cannot be of finite order. We consider the other two possibilities.

\textbf{Case 1}\qua $f$ is pseudo-Anosov.

Let $\!\lai f \rai\!$ denote the cyclic subgroup of $\mathcal{A}$ generated by $f$, and let $\mathcal{C}\!=\!C_{\mathcal{M}(\mathbf{S})}(\lai f \rai )$, the centralizer of $\lai f \rai $ in $\mathcal{M}(\mathbf{S})$.  Then $\mathcal{A} \subseteq \mathcal{C}$ and $\mathcal{A}$ is torsion-free. We conclude by a theorem of McCarthy (\cite{4}, Corollary~3) that $\mathcal{A}$ is infinite cyclic.  Hence $\mathrm{rank}(\mathcal{A}) = 1 \leq 2g-3$.

\textbf{Case 2}\qua $f$ is reducible.

Given $h \in \mathcal{A}$, let $\fr{E}_h$ denote the essential reduction system of $h$, and let \[ \fr{E}= \bigcup_{h \in \mathcal {A}} \fr{E}_h \]
Then $\fr{E}$ is an adequate reduction system for $\mathcal{A}$ (\cite{1}, Lemma~3.1(1)), and f reducible implies $\fr{E} \neq \emptyset$, so every element of $\mathcal{A}$ is reducible.

Let $\mathcal{M}_\fr{E}(\mathbf{S})$ denote the stabilizer of $\fr{E}$ in $\mathcal{M}(\mathbf{S})$, and let $\Lambda \co \mathcal{M}_\fr{E}(\mathbf{S}) \to \mathcal{M}(\mathbf{S}_\fr{E})$ be the reduction homomorphism.  Then $\mathrm{ker}(\Lambda) = \mathcal{D}_\fr{E}$, the multitwist group on $\fr{E}$, and thus \[ \mathrm{ker}(\Lambda \arrowvert_\mathcal{A}) = \mathrm{ker}(\Lambda) \cap \mathcal{A} = \mathcal{D}_\fr{E} \cap \mathcal{A}=\mathcal{D}_\fr{E} \cap \mathcal{T} \cap \mathcal{A} =\mathcal{T}_\fr{E} \cap \mathcal{A}. \]
We now have a short exact sequence \[ 0 \longrightarrow \mathcal{T}_\fr{E} \cap \mathcal{A} \longrightarrow \mathcal{A} \stackrel{\Lambda \arrowvert_{\mathcal{A}}}{\longrightarrow}\Lambda(\mathcal{A}) \longrightarrow 0 \] of free Abelian groups, which shows that \[ \mathrm{rank}(\mathcal{A}) = \mathrm{rank}(\mathcal{T}_\fr{E} \cap \mathcal{A}) + \mathrm{rank}(\Lambda(\mathcal{A})) \leq \mathrm{rank}(\mathcal{T}_\fr{E}) + \mathrm{rank}(\Lambda(\mathcal{A})). \]
We will be done, by applying Lemma~3.2, once we show that $\mathrm{rank}(\Lambda(\mathcal{A})) \leq \Omega$, the number of components of $\scut$ not homeomorphic to a pair of pants or a one-holed torus.

A theorem of Ivanov (\cite{3}, Theorem~1.2) implies that $\Lambda(f)$ restricts to each component $\mathbf{S}_1,\mathbf{S}_2, \dots ,\mathbf{S}_\nu$ of $\scut$, giving ``projections'' $ p_i \co  \Lambda(\mathcal{A}) \longrightarrow\mathcal{M}(\mathbf{S}_i)$ induced by restricting representatives.  Set $\mathcal{A}_i=p_i(\Lambda(\mathcal{A})) \subseteq \mathcal{M}(\mathbf{S}_i)$.  Then $\Lambda(\mathcal{A}) \subseteq \bigoplus \mathcal{A}_i$, so $\mathrm{rank}(\Lambda(\mathcal{A})) \leq \sum \mathrm{rank}(\mathcal{A}_i)$.  We make the following observations:

(i)\qua  If $\mathbf{S}_i$ is a pair of pants, then $\mathcal{M}(\mathbf{S}_i)$ is finite, so $\mathrm{rank}(\mathcal{A}_i) = 0$.

(ii)\qua  If $\mathbf{S}_i$ is a one-holed torus, then the homomorphism $H_1(\mathbf{S}_i) \to H_1(\mathbf{S})$ induced by inclusion is injective. Any homeomorphism $f$ representing an element of $\mathcal{A}$ maps a circle $\fr{c}$ in $\mathbf{S}_i$ to a circle $\fr{c'}$ in $\mathbf{S}_i$, so $\mathcal{A}_i$ lies within the Torelli group of $\mathbf{S}_i$, which is trivial in this case.

(iii)\qua  If $\mathbf{S}_i$ is neither a pair of pants nor a one-holed torus, then $\mathcal{A}_i$ is either trivial or is an adequately reduced torsion-free Abelian subgroup of $\mathcal{M}(\mathbf{S}_i)$. So again by McCarthy's theorem, $\mathrm{rank}(\mathcal{A}_i) \leq 1$.

These observations tell us that
$$\mathrm{rank}(\Lambda(\mathcal{A})) \leq \sum_{i=1}^\nu \mathrm{rank}(\mathcal{A}_i) \leq \Omega.\eqno{\qed}$$

\Addresses\recd

\end{document}

%% file: agtout.tex
\def\ifplaintex{\expandafter\ifx\csname documentclass\endcsname\relax}

\def\gtp{{\mathsurround=0pt\it $\cal G\mskip-2mu$eometry \&\ 
$\cal T\!\!$opology $\cal P\!$ublications}}  

\def\recd{{\small Received:\qua\receiveddate\ifx\reviseddate\relax
\else\qquad Revised:\qua\reviseddate\fi\par}} 

\def\lognumber#1{\def\thelognumber{#1}}
\def\volumenumber#1{\def\thevolumenumber{#1}}
\def\volumeyear#1{\def\thevolumeyear{#1}}
\def\papernumber#1{\def\thepapernumber{#1}}
\def\pagenumbers#1#2{\def\startpage{#1}\def\finishpage{#2}}
\def\published#1{\def\publishdate{#1}}

\def\received#1{\def\receiveddate{#1}}
\def\revised#1{\def\reviseddate{#1}}
\def\accepted#1{\def\accepteddate{#1}}

\long\def\asciiabstract#1{\long\def\theasciiabstract{#1}}

\let\\\par\let\thelognumber\relax\let\thevolumenumber\relax
\let\thepapernumber\relax\let\thevolumeyear\relax\let\startpage\relax
\let\finishpage\relax\let\publishdate\relax\let\receiveddate\relax
\let\reviseddate\relax\let\accepteddate\relax\let\theasciititle\relax
\let\theasciiauthors\relax
\let\theasciiabstract\relax

\let\theasciiemail\relax

\ifplaintex
\font\logobig=cmssbx10 scaled 3836
\font\logomed=cmssbx10 scaled 2557
\else
\font\logobig=cmssbx10 scaled 4200
\font\logomed=cmssbx10 scaled 2800
\fi

\long\def\makeagttitle{   
\count0=\startpage
\agt\hfill      
\hbox to 45truept{\vbox to 0pt{\vglue -13truept{\logomed A\kern -.37em{\logobig 
T}\kern -.38em G}\vss}\hss}
\break
{\small Volume \thevolumenumber\ (\thevolumeyear)
\startpage--\finishpage\nl
Published: \publishdate}

\vglue .25truein

{\parskip=0pt\leftskip 0pt plus
1fil\def\\{\par\smallskip}{\Large\bf\thetitle}\par\medskip} \vglue
0.05truein

{\parskip=0pt\leftskip 0pt plus 1fil\def\\{\par}{\sc\theauthors}
\par\medskip}%
 
\vglue 0.03truein 

{\small\leftskip 25truept\rightskip 25truept{\bf Abstract}\stdspace\theabstract

{\bf AMS Classification}\stdspace\theprimaryclass
\ifx\thesecondaryclass\relax\else; \thesecondaryclass\fi\par
{\bf Keywords}\stdspace \thekeywords\par}\vglue 7truept

}   

\ifplaintex
\font\phead=cmsl9 scaled 950
\font\pnum=cmbx10 scaled 913
\font\pfoot=cmsl9 scaled 950
\headline{\vbox to 0pt{\vskip -4.5mm\line{\small\phead\ifnum
\count0=\startpage ISSN 1472-2739 (on-line) 1472-2747 (printed)
\hfill {\pnum\folio}\else\ifodd\count0\def\\{ }%
\ifx\theshorttitle\relax\thetitle\else\theshorttitle\fi\hfill{\pnum\folio}
\else\def\\{ and }{\pnum\folio}\hfill\ifx\theshortauthors\relax\theauthors
\else\theshortauthors\fi\fi\fi}\vss}}
\footline{\vbox to 0pt{\vglue 0mm\line{\small\pfoot\ifnum\count0=\startpage
\copyright\ \gtp\hfill\else
\agt, Volume \thevolumenumber\ (\thevolumeyear)\hfill\fi}\vss}}
\else
\font=cmsl9 scaled 1050
\font\lnum=cmbx10 
\font=cmsl9 scaled 1050
\makeatletter
\def\@oddhead{{\small\lhead\ifnum\count0=\startpage ISSN 1472-2739 
(on-line) 1472-2747 (printed)\hfill {\lnum\number\count0}\else\ifodd\count0
\def\\{ }\ifx\theshorttitle\relax \thetitle \else\theshorttitle\fi\hfill
{\lnum\number\count0}\else\def\\{ and }{\lnum\number\count0}
\hfill\ifx\theshortauthors\relax 
\theauthors\else\theshortauthors\fi\fi\fi}}\def\@evenhead{\@oddhead}
\def\@oddfoot{\small\lfoot\ifnum\count0=\startpage\copyright\ \gtp\hfill\else
\agt, Volume \thevolumenumber\ (\thevolumeyear)\hfill\fi}
\def\@evenfoot{\@oddfoot}
\makeatother
\fi
\let\maketitlepage\makeagttitle
\let\makeshorttitle\maketitlepage
\let\maketitle\maketitlepage

\newwrite\gtoutfile
\long\gdef\makeheadfile{  
{\def\\{, }\def\s{ }
\immediate\openout\gtoutfile head.xxx
\immediate\write\gtoutfile{To: math@arxiv.org}
\immediate\write\gtoutfile{Subject: put OR rep NNNNN:ppppp}
\immediate\write\gtoutfile{--text follows this line--}
\immediate\write\gtoutfile{Proxy-for: \ifx\theasciiauthors\relax
\theauthors\else\theasciiauthors\fi\s<\ifx\theasciiemail\relax\theemail\else\theasciiemail\fi>}
\immediate\write\gtoutfile{\noexpand\\}
\immediate\write\gtoutfile{Authors: \ifx\theasciiauthors\relax
\theauthors\else\theasciiauthors\fi}
{\def\\{ }\immediate\write\gtoutfile{Title: \ifx\theasciititle\relax
\thetitle\else\theasciititle\fi}}
\immediate\write\gtoutfile{Subj-class: GT or SG, GR etc}
\immediate\write\gtoutfile{MSC-class: \theprimaryclass\ifx\thesecondaryclass\relax\else, \thesecondaryclass\fi}
\immediate\write\gtoutfile{Journal-ref: Algebr. Geom. Topol. \thevolumenumber\s
(\thevolumeyear) \startpage-\finishpage}
\immediate\write\gtoutfile{Comments: Published by Algebraic and
Geometric Topology at}
\immediate\write\gtoutfile{\s\s\s  http://www.maths.warwick.ac.uk/agt/AGTVol\thevolumenumber/agt-\thevolumenumber-\thepapernumber.abs.html}
\immediate\write\gtoutfile{\noexpand\\}
\immediate\write\gtoutfile{}
\ifx\theasciiabstract\relax
\immediate\write\gtoutfile{\theabstract}\else
\immediate\write\gtoutfile{\theasciiabstract}\fi
\immediate\write\gtoutfile{}
\immediate\write\gtoutfile{\noexpand\\}
\immediate\write\gtoutfile{}
\immediate\closeout\gtoutfile}}  

\def\maketitlepage{\makeagttitle\makeheadfile}
\let\makeshorttitle\maketitlepage
\let\maketitle\maketitlepage

\def\ifplaintex{\expandafter\ifx\csname documentclass\endcsname\relax}

\def\gtp{{\mathsurround=0pt\it $\cal G\mskip-2mu$eometry \&\ 
$\cal T\!\!$opology $\cal P\!$ublications}}  

\def\recd{{\small Received:\qua\receiveddate\ifx\reviseddate\relax
\else\qquad Revised:\qua\reviseddate\fi\par}} 

\def\lognumber#1{\def\thelognumber{#1}}
\def\volumenumber#1{\def\thevolumenumber{#1}}
\def\volumeyear#1{\def\thevolumeyear{#1}}
\def\papernumber#1{\def\thepapernumber{#1}}
\def\pagenumbers#1#2{\def\startpage{#1}\def\finishpage{#2}}
\def\published#1{\def\publishdate{#1}}

\def\received#1{\def\receiveddate{#1}}
\def\revised#1{\def\reviseddate{#1}}
\def\accepted#1{\def\accepteddate{#1}}

\long\def\asciiabstract#1{\long\def\theasciiabstract{#1}}

\let\\\par\let\thelognumber\relax\let\thevolumenumber\relax
\let\thepapernumber\relax\let\thevolumeyear\relax\let\startpage\relax
\let\finishpage\relax\let\publishdate\relax\let\receiveddate\relax
\let\reviseddate\relax\let\accepteddate\relax\let\theasciititle\relax
\let\theasciiauthors\relax
\let\theasciiabstract\relax

\let\theasciiemail\relax

\ifplaintex
\font\logobig=cmssbx10 scaled 3836
\font\logomed=cmssbx10 scaled 2557
\else
\font\logobig=cmssbx10 scaled 4200
\font\logomed=cmssbx10 scaled 2800
\fi

\long\def\makeagttitle{   
\count0=\startpage
\agt\hfill      
\hbox to 45truept{\vbox to 0pt{\vglue -13truept{\logomed A\kern -.37em{\logobig 
T}\kern -.38em G}\vss}\hss}
\break
{\small Volume \thevolumenumber\ (\thevolumeyear)
\startpage--\finishpage\nl
Published: \publishdate}

\vglue .25truein

{\parskip=0pt\leftskip 0pt plus
1fil\def\\{\par\smallskip}{\Large\bf\thetitle}\par\medskip} \vglue
0.05truein

{\parskip=0pt\leftskip 0pt plus 1fil\def\\{\par}{\sc\theauthors}
\par\medskip}%
 
\vglue 0.03truein 

{\small\leftskip 25truept\rightskip 25truept{\bf Abstract}\stdspace\theabstract

{\bf AMS Classification}\stdspace\theprimaryclass
\ifx\thesecondaryclass\relax\else; \thesecondaryclass\fi\par
{\bf Keywords}\stdspace \thekeywords\par}\vglue 7truept

}   

\ifplaintex
\font\phead=cmsl9 scaled 950
\font\pnum=cmbx10 scaled 913
\font\pfoot=cmsl9 scaled 950
\headline{\vbox to 0pt{\vskip -4.5mm\line{\small\phead\ifnum
\count0=\startpage ISSN 1472-2739 (on-line) 1472-2747 (printed)
\hfill {\pnum\folio}\else\ifodd\count0\def\\{ }%
\ifx\theshorttitle\relax\thetitle\else\theshorttitle\fi\hfill{\pnum\folio}
\else\def\\{ and }{\pnum\folio}\hfill\ifx\theshortauthors\relax\theauthors
\else\theshortauthors\fi\fi\fi}\vss}}
\footline{\vbox to 0pt{\vglue 0mm\line{\small\pfoot\ifnum\count0=\startpage
\copyright\ \gtp\hfill\else
\agt, Volume \thevolumenumber\ (\thevolumeyear)\hfill\fi}\vss}}
\else
\font=cmsl9 scaled 1050
\font\lnum=cmbx10 
\font=cmsl9 scaled 1050
\makeatletter
\def\@oddhead{{\small\lhead\ifnum\count0=\startpage ISSN 1472-2739 
(on-line) 1472-2747 (printed)\hfill {\lnum\number\count0}\else\ifodd\count0
\def\\{ }\ifx\theshorttitle\relax \thetitle \else\theshorttitle\fi\hfill
{\lnum\number\count0}\else\def\\{ and }{\lnum\number\count0}
\hfill\ifx\theshortauthors\relax 
\theauthors\else\theshortauthors\fi\fi\fi}}\def\@evenhead{\@oddhead}
\def\@oddfoot{\small\lfoot\ifnum\count0=\startpage\copyright\ \gtp\hfill\else
\agt, Volume \thevolumenumber\ (\thevolumeyear)\hfill\fi}
\def\@evenfoot{\@oddfoot}
\makeatother
\fi
\let\maketitlepage\makeagttitle
\let\makeshorttitle\maketitlepage
\let\maketitle\maketitlepage

\newwrite\gtoutfile
\long\gdef\makeheadfile{  
{\def\\{, }\def\s{ }
\immediate\openout\gtoutfile head.xxx
\immediate\write\gtoutfile{To: math@arxiv.org}
\immediate\write\gtoutfile{Subject: put OR rep NNNNN:ppppp}
\immediate\write\gtoutfile{--text follows this line--}
\immediate\write\gtoutfile{Proxy-for: \ifx\theasciiauthors\relax
\theauthors\else\theasciiauthors\fi\s<\ifx\theasciiemail\relax\theemail\else\theasciiemail\fi>}
\immediate\write\gtoutfile{\noexpand\\}
\immediate\write\gtoutfile{Authors: \ifx\theasciiauthors\relax
\theauthors\else\theasciiauthors\fi}
{\def\\{ }\immediate\write\gtoutfile{Title: \ifx\theasciititle\relax
\thetitle\else\theasciititle\fi}}
\immediate\write\gtoutfile{Subj-class: GT or SG, GR etc}
\immediate\write\gtoutfile{MSC-class: \theprimaryclass\ifx\thesecondaryclass\relax\else, \thesecondaryclass\fi}
\immediate\write\gtoutfile{Journal-ref: Algebr. Geom. Topol. \thevolumenumber\s
(\thevolumeyear) \startpage-\finishpage}
\immediate\write\gtoutfile{Comments: Published by Algebraic and
Geometric Topology at}
\immediate\write\gtoutfile{\s\s\s  http://www.maths.warwick.ac.uk/agt/AGTVol\thevolumenumber/agt-\thevolumenumber-\thepapernumber.abs.html}
\immediate\write\gtoutfile{\noexpand\\}
\immediate\write\gtoutfile{}
\ifx\theasciiabstract\relax
\immediate\write\gtoutfile{\theabstract}\else
\immediate\write\gtoutfile{\theasciiabstract}\fi
\immediate\write\gtoutfile{}
\immediate\write\gtoutfile{\noexpand\\}
\immediate\write\gtoutfile{}
\immediate\closeout\gtoutfile}}  

\def\maketitlepage{\makeagttitle\makeheadfile}
\let\makeshorttitle\maketitlepage
\let\maketitle\maketitlepage

\def\ifplaintex{\expandafter\ifx\csname documentclass\endcsname\relax}

\def\gtp{{\mathsurround=0pt\it $\cal G\mskip-2mu$eometry \&\ 
$\cal T\!\!$opology $\cal P\!$ublications}}  

\def\recd{{\small Received:\qua\receiveddate\ifx\reviseddate\relax
\else\qquad Revised:\qua\reviseddate\fi\par}} 

\def\lognumber#1{\def\thelognumber{#1}}
\def\volumenumber#1{\def\thevolumenumber{#1}}
\def\volumeyear#1{\def\thevolumeyear{#1}}
\def\papernumber#1{\def\thepapernumber{#1}}
\def\pagenumbers#1#2{\def\startpage{#1}\def\finishpage{#2}}
\def\published#1{\def\publishdate{#1}}

\def\received#1{\def\receiveddate{#1}}
\def\revised#1{\def\reviseddate{#1}}
\def\accepted#1{\def\accepteddate{#1}}

\long\def\asciiabstract#1{\long\def\theasciiabstract{#1}}

\let\\\par\let\thelognumber\relax\let\thevolumenumber\relax
\let\thepapernumber\relax\let\thevolumeyear\relax\let\startpage\relax
\let\finishpage\relax\let\publishdate\relax\let\receiveddate\relax
\let\reviseddate\relax\let\accepteddate\relax\let\theasciititle\relax
\let\theasciiauthors\relax
\let\theasciiabstract\relax

\let\theasciiemail\relax

\ifplaintex
\font\logobig=cmssbx10 scaled 3836
\font\logomed=cmssbx10 scaled 2557
\else
\font\logobig=cmssbx10 scaled 4200
\font\logomed=cmssbx10 scaled 2800
\fi

\long\def\makeagttitle{   
\count0=\startpage
\agt\hfill      
\hbox to 45truept{\vbox to 0pt{\vglue -13truept{\logomed A\kern -.37em{\logobig 
T}\kern -.38em G}\vss}\hss}
\break
{\small Volume \thevolumenumber\ (\thevolumeyear)
\startpage--\finishpage\nl
Published: \publishdate}

\vglue .25truein

{\parskip=0pt\leftskip 0pt plus
1fil\def\\{\par\smallskip}{\Large\bf\thetitle}\par\medskip} \vglue
0.05truein

{\parskip=0pt\leftskip 0pt plus 1fil\def\\{\par}{\sc\theauthors}
\par\medskip}%
 
\vglue 0.03truein 

{\small\leftskip 25truept\rightskip 25truept{\bf Abstract}\stdspace\theabstract

{\bf AMS Classification}\stdspace\theprimaryclass
\ifx\thesecondaryclass\relax\else; \thesecondaryclass\fi\par
{\bf Keywords}\stdspace \thekeywords\par}\vglue 7truept

}   

\ifplaintex
\font\phead=cmsl9 scaled 950
\font\pnum=cmbx10 scaled 913
\font\pfoot=cmsl9 scaled 950
\headline{\vbox to 0pt{\vskip -4.5mm\line{\small\phead\ifnum
\count0=\startpage ISSN 1472-2739 (on-line) 1472-2747 (printed)
\hfill {\pnum\folio}\else\ifodd\count0\def\\{ }%
\ifx\theshorttitle\relax\thetitle\else\theshorttitle\fi\hfill{\pnum\folio}
\else\def\\{ and }{\pnum\folio}\hfill\ifx\theshortauthors\relax\theauthors
\else\theshortauthors\fi\fi\fi}\vss}}
\footline{\vbox to 0pt{\vglue 0mm\line{\small\pfoot\ifnum\count0=\startpage
\copyright\ \gtp\hfill\else
\agt, Volume \thevolumenumber\ (\thevolumeyear)\hfill\fi}\vss}}
\else
\font=cmsl9 scaled 1050
\font\lnum=cmbx10 
\font=cmsl9 scaled 1050
\makeatletter
\def\@oddhead{{\small\lhead\ifnum\count0=\startpage ISSN 1472-2739 
(on-line) 1472-2747 (printed)\hfill {\lnum\number\count0}\else\ifodd\count0
\def\\{ }\ifx\theshorttitle\relax \thetitle \else\theshorttitle\fi\hfill
{\lnum\number\count0}\else\def\\{ and }{\lnum\number\count0}
\hfill\ifx\theshortauthors\relax 
\theauthors\else\theshortauthors\fi\fi\fi}}\def\@evenhead{\@oddhead}
\def\@oddfoot{\small\lfoot\ifnum\count0=\startpage\copyright\ \gtp\hfill\else
\agt, Volume \thevolumenumber\ (\thevolumeyear)\hfill\fi}
\def\@evenfoot{\@oddfoot}
\makeatother
\fi
\let\maketitlepage\makeagttitle
\let\makeshorttitle\maketitlepage
\let\maketitle\maketitlepage

\newwrite\gtoutfile
\long\gdef\makeheadfile{  
{\def\\{, }\def\s{ }
\immediate\openout\gtoutfile head.xxx
\immediate\write\gtoutfile{To: math@arxiv.org}
\immediate\write\gtoutfile{Subject: put OR rep NNNNN:ppppp}
\immediate\write\gtoutfile{--text follows this line--}
\immediate\write\gtoutfile{Proxy-for: \ifx\theasciiauthors\relax
\theauthors\else\theasciiauthors\fi\s<\ifx\theasciiemail\relax\theemail\else\theasciiemail\fi>}
\immediate\write\gtoutfile{\noexpand\\}
\immediate\write\gtoutfile{Authors: \ifx\theasciiauthors\relax
\theauthors\else\theasciiauthors\fi}
{\def\\{ }\immediate\write\gtoutfile{Title: \ifx\theasciititle\relax
\thetitle\else\theasciititle\fi}}
\immediate\write\gtoutfile{Subj-class: GT or SG, GR etc}
\immediate\write\gtoutfile{MSC-class: \theprimaryclass\ifx\thesecondaryclass\relax\else, \thesecondaryclass\fi}
\immediate\write\gtoutfile{Journal-ref: Algebr. Geom. Topol. \thevolumenumber\s
(\thevolumeyear) \startpage-\finishpage}
\immediate\write\gtoutfile{Comments: Published by Algebraic and
Geometric Topology at}
\immediate\write\gtoutfile{\s\s\s  http://www.maths.warwick.ac.uk/agt/AGTVol\thevolumenumber/agt-\thevolumenumber-\thepapernumber.abs.html}
\immediate\write\gtoutfile{\noexpand\\}
\immediate\write\gtoutfile{}
\ifx\theasciiabstract\relax
\immediate\write\gtoutfile{\theabstract}\else
\immediate\write\gtoutfile{\theasciiabstract}\fi
\immediate\write\gtoutfile{}
\immediate\write\gtoutfile{\noexpand\\}
\immediate\write\gtoutfile{}
\immediate\closeout\gtoutfile}}  

\def\maketitlepage{\makeagttitle\makeheadfile}
\let\makeshorttitle\maketitlepage
\let\maketitle\maketitlepage

\def\ifplaintex{\expandafter\ifx\csname documentclass\endcsname\relax}

\def\gtp{{\mathsurround=0pt\it $\cal G\mskip-2mu$eometry \&\ 
$\cal T\!\!$opology $\cal P\!$ublications}}  

\def\recd{{\small Received:\qua\receiveddate\ifx\reviseddate\relax
\else\qquad Revised:\qua\reviseddate\fi\par}} 

\def\lognumber#1{\def\thelognumber{#1}}
\def\volumenumber#1{\def\thevolumenumber{#1}}
\def\volumeyear#1{\def\thevolumeyear{#1}}
\def\papernumber#1{\def\thepapernumber{#1}}
\def\pagenumbers#1#2{\def\startpage{#1}\def\finishpage{#2}}
\def\published#1{\def\publishdate{#1}}

\def\received#1{\def\receiveddate{#1}}
\def\revised#1{\def\reviseddate{#1}}
\def\accepted#1{\def\accepteddate{#1}}

\long\def\asciiabstract#1{\long\def\theasciiabstract{#1}}

\let\\\par\let\thelognumber\relax\let\thevolumenumber\relax
\let\thepapernumber\relax\let\thevolumeyear\relax\let\startpage\relax
\let\finishpage\relax\let\publishdate\relax\let\receiveddate\relax
\let\reviseddate\relax\let\accepteddate\relax\let\theasciititle\relax
\let\theasciiauthors\relax
\let\theasciiabstract\relax

\let\theasciiemail\relax

\ifplaintex
\font\logobig=cmssbx10 scaled 3836
\font\logomed=cmssbx10 scaled 2557
\else
\font\logobig=cmssbx10 scaled 4200
\font\logomed=cmssbx10 scaled 2800
\fi

\long\def\makeagttitle{   
\count0=\startpage
\agt\hfill      
\hbox to 45truept{\vbox to 0pt{\vglue -13truept{\logomed A\kern -.37em{\logobig 
T}\kern -.38em G}\vss}\hss}
\break
{\small Volume \thevolumenumber\ (\thevolumeyear)
\startpage--\finishpage\nl
Published: \publishdate}

\vglue .25truein

{\parskip=0pt\leftskip 0pt plus
1fil\def\\{\par\smallskip}{\Large\bf\thetitle}\par\medskip} \vglue
0.05truein

{\parskip=0pt\leftskip 0pt plus 1fil\def\\{\par}{\sc\theauthors}
\par\medskip}%
 
\vglue 0.03truein 

{\small\leftskip 25truept\rightskip 25truept{\bf Abstract}\stdspace\theabstract

{\bf AMS Classification}\stdspace\theprimaryclass
\ifx\thesecondaryclass\relax\else; \thesecondaryclass\fi\par
{\bf Keywords}\stdspace \thekeywords\par}\vglue 7truept

}   

\ifplaintex
\font\phead=cmsl9 scaled 950
\font\pnum=cmbx10 scaled 913
\font\pfoot=cmsl9 scaled 950
\headline{\vbox to 0pt{\vskip -4.5mm\line{\small\phead\ifnum
\count0=\startpage ISSN 1472-2739 (on-line) 1472-2747 (printed)
\hfill {\pnum\folio}\else\ifodd\count0\def\\{ }%
\ifx\theshorttitle\relax\thetitle\else\theshorttitle\fi\hfill{\pnum\folio}
\else\def\\{ and }{\pnum\folio}\hfill\ifx\theshortauthors\relax\theauthors
\else\theshortauthors\fi\fi\fi}\vss}}
\footline{\vbox to 0pt{\vglue 0mm\line{\small\pfoot\ifnum\count0=\startpage
\copyright\ \gtp\hfill\else
\agt, Volume \thevolumenumber\ (\thevolumeyear)\hfill\fi}\vss}}
\else
\font=cmsl9 scaled 1050
\font\lnum=cmbx10 
\font=cmsl9 scaled 1050
\makeatletter
\def\@oddhead{{\small\lhead\ifnum\count0=\startpage ISSN 1472-2739 
(on-line) 1472-2747 (printed)\hfill {\lnum\number\count0}\else\ifodd\count0
\def\\{ }\ifx\theshorttitle\relax \thetitle \else\theshorttitle\fi\hfill
{\lnum\number\count0}\else\def\\{ and }{\lnum\number\count0}
\hfill\ifx\theshortauthors\relax 
\theauthors\else\theshortauthors\fi\fi\fi}}\def\@evenhead{\@oddhead}
\def\@oddfoot{\small\lfoot\ifnum\count0=\startpage\copyright\ \gtp\hfill\else
\agt, Volume \thevolumenumber\ (\thevolumeyear)\hfill\fi}
\def\@evenfoot{\@oddfoot}
\makeatother
\fi
\let\maketitlepage\makeagttitle
\let\makeshorttitle\maketitlepage
\let\maketitle\maketitlepage

\newwrite\gtoutfile
\long\gdef\makeheadfile{  
{\def\\{, }\def\s{ }
\immediate\openout\gtoutfile head.xxx
\immediate\write\gtoutfile{To: math@arxiv.org}
\immediate\write\gtoutfile{Subject: put OR rep NNNNN:ppppp}
\immediate\write\gtoutfile{--text follows this line--}
\immediate\write\gtoutfile{Proxy-for: \ifx\theasciiauthors\relax
\theauthors\else\theasciiauthors\fi\s<\ifx\theasciiemail\relax\theemail\else\theasciiemail\fi>}
\immediate\write\gtoutfile{\noexpand\\}
\immediate\write\gtoutfile{Authors: \ifx\theasciiauthors\relax
\theauthors\else\theasciiauthors\fi}
{\def\\{ }\immediate\write\gtoutfile{Title: \ifx\theasciititle\relax
\thetitle\else\theasciititle\fi}}
\immediate\write\gtoutfile{Subj-class: GT or SG, GR etc}
\immediate\write\gtoutfile{MSC-class: \theprimaryclass\ifx\thesecondaryclass\relax\else, \thesecondaryclass\fi}
\immediate\write\gtoutfile{Journal-ref: Algebr. Geom. Topol. \thevolumenumber\s
(\thevolumeyear) \startpage-\finishpage}
\immediate\write\gtoutfile{Comments: Published by Algebraic and
Geometric Topology at}
\immediate\write\gtoutfile{\s\s\s  http://www.maths.warwick.ac.uk/agt/AGTVol\thevolumenumber/agt-\thevolumenumber-\thepapernumber.abs.html}
\immediate\write\gtoutfile{\noexpand\\}
\immediate\write\gtoutfile{}
\ifx\theasciiabstract\relax
\immediate\write\gtoutfile{\theabstract}\else
\immediate\write\gtoutfile{\theasciiabstract}\fi
\immediate\write\gtoutfile{}
\immediate\write\gtoutfile{\noexpand\\}
\immediate\write\gtoutfile{}
\immediate\closeout\gtoutfile}}  

\def\maketitlepage{\makeagttitle\makeheadfile}
\let\makeshorttitle\maketitlepage
\let\maketitle\maketitlepage

\def\ifplaintex{\expandafter\ifx\csname documentclass\endcsname\relax}

\def\gtp{{\mathsurround=0pt\it $\cal G\mskip-2mu$eometry \&\ 
$\cal T\!\!$opology $\cal P\!$ublications}}  

\def\recd{{\small Received:\qua\receiveddate\ifx\reviseddate\relax
\else\qquad Revised:\qua\reviseddate\fi\par}} 

\def\lognumber#1{\def\thelognumber{#1}}
\def\volumenumber#1{\def\thevolumenumber{#1}}
\def\volumeyear#1{\def\thevolumeyear{#1}}
\def\papernumber#1{\def\thepapernumber{#1}}
\def\pagenumbers#1#2{\def\startpage{#1}\def\finishpage{#2}}
\def\published#1{\def\publishdate{#1}}

\def\received#1{\def\receiveddate{#1}}
\def\revised#1{\def\reviseddate{#1}}
\def\accepted#1{\def\accepteddate{#1}}

\long\def\asciiabstract#1{\long\def\theasciiabstract{#1}}

\let\\\par\let\thelognumber\relax\let\thevolumenumber\relax
\let\thepapernumber\relax\let\thevolumeyear\relax\let\startpage\relax
\let\finishpage\relax\let\publishdate\relax\let\receiveddate\relax
\let\reviseddate\relax\let\accepteddate\relax\let\theasciititle\relax
\let\theasciiauthors\relax
\let\theasciiabstract\relax

\let\theasciiemail\relax

\ifplaintex
\font\logobig=cmssbx10 scaled 3836
\font\logomed=cmssbx10 scaled 2557
\else
\font\logobig=cmssbx10 scaled 4200
\font\logomed=cmssbx10 scaled 2800
\fi

\long\def\makeagttitle{   
\count0=\startpage
\agt\hfill      
\hbox to 45truept{\vbox to 0pt{\vglue -13truept{\logomed A\kern -.37em{\logobig 
T}\kern -.38em G}\vss}\hss}
\break
{\small Volume \thevolumenumber\ (\thevolumeyear)
\startpage--\finishpage\nl
Published: \publishdate}

\vglue .25truein

{\parskip=0pt\leftskip 0pt plus
1fil\def\\{\par\smallskip}{\Large\bf\thetitle}\par\medskip} \vglue
0.05truein

{\parskip=0pt\leftskip 0pt plus 1fil\def\\{\par}{\sc\theauthors}
\par\medskip}%
 
\vglue 0.03truein 

{\small\leftskip 25truept\rightskip 25truept{\bf Abstract}\stdspace\theabstract

{\bf AMS Classification}\stdspace\theprimaryclass
\ifx\thesecondaryclass\relax\else; \thesecondaryclass\fi\par
{\bf Keywords}\stdspace \thekeywords\par}\vglue 7truept

}   

\ifplaintex
\font\phead=cmsl9 scaled 950
\font\pnum=cmbx10 scaled 913
\font\pfoot=cmsl9 scaled 950
\headline{\vbox to 0pt{\vskip -4.5mm\line{\small\phead\ifnum
\count0=\startpage ISSN 1472-2739 (on-line) 1472-2747 (printed)
\hfill {\pnum\folio}\else\ifodd\count0\def\\{ }%
\ifx\theshorttitle\relax\thetitle\else\theshorttitle\fi\hfill{\pnum\folio}
\else\def\\{ and }{\pnum\folio}\hfill\ifx\theshortauthors\relax\theauthors
\else\theshortauthors\fi\fi\fi}\vss}}
\footline{\vbox to 0pt{\vglue 0mm\line{\small\pfoot\ifnum\count0=\startpage
\copyright\ \gtp\hfill\else
\agt, Volume \thevolumenumber\ (\thevolumeyear)\hfill\fi}\vss}}
\else
\font=cmsl9 scaled 1050
\font\lnum=cmbx10 
\font=cmsl9 scaled 1050
\makeatletter
\def\@oddhead{{\small\lhead\ifnum\count0=\startpage ISSN 1472-2739 
(on-line) 1472-2747 (printed)\hfill {\lnum\number\count0}\else\ifodd\count0
\def\\{ }\ifx\theshorttitle\relax \thetitle \else\theshorttitle\fi\hfill
{\lnum\number\count0}\else\def\\{ and }{\lnum\number\count0}
\hfill\ifx\theshortauthors\relax 
\theauthors\else\theshortauthors\fi\fi\fi}}\def\@evenhead{\@oddhead}
\def\@oddfoot{\small\lfoot\ifnum\count0=\startpage\copyright\ \gtp\hfill\else
\agt, Volume \thevolumenumber\ (\thevolumeyear)\hfill\fi}
\def\@evenfoot{\@oddfoot}
\makeatother
\fi
\let\maketitlepage\makeagttitle
\let\makeshorttitle\maketitlepage
\let\maketitle\maketitlepage

\newwrite\gtoutfile
\long\gdef\makeheadfile{  
{\def\\{, }\def\s{ }
\immediate\openout\gtoutfile head.xxx
\immediate\write\gtoutfile{To: math@arxiv.org}
\immediate\write\gtoutfile{Subject: put OR rep NNNNN:ppppp}
\immediate\write\gtoutfile{--text follows this line--}
\immediate\write\gtoutfile{Proxy-for: \ifx\theasciiauthors\relax
\theauthors\else\theasciiauthors\fi\s<\ifx\theasciiemail\relax\theemail\else\theasciiemail\fi>}
\immediate\write\gtoutfile{\noexpand\\}
\immediate\write\gtoutfile{Authors: \ifx\theasciiauthors\relax
\theauthors\else\theasciiauthors\fi}
{\def\\{ }\immediate\write\gtoutfile{Title: \ifx\theasciititle\relax
\thetitle\else\theasciititle\fi}}
\immediate\write\gtoutfile{Subj-class: GT or SG, GR etc}
\immediate\write\gtoutfile{MSC-class: \theprimaryclass\ifx\thesecondaryclass\relax\else, \thesecondaryclass\fi}
\immediate\write\gtoutfile{Journal-ref: Algebr. Geom. Topol. \thevolumenumber\s
(\thevolumeyear) \startpage-\finishpage}
\immediate\write\gtoutfile{Comments: Published by Algebraic and
Geometric Topology at}
\immediate\write\gtoutfile{\s\s\s  http://www.maths.warwick.ac.uk/agt/AGTVol\thevolumenumber/agt-\thevolumenumber-\thepapernumber.abs.html}
\immediate\write\gtoutfile{\noexpand\\}
\immediate\write\gtoutfile{}
\ifx\theasciiabstract\relax
\immediate\write\gtoutfile{\theabstract}\else
\immediate\write\gtoutfile{\theasciiabstract}\fi
\immediate\write\gtoutfile{}
\immediate\write\gtoutfile{\noexpand\\}
\immediate\write\gtoutfile{}
\immediate\closeout\gtoutfile}}  

\def\maketitlepage{\makeagttitle\makeheadfile}
\let\makeshorttitle\maketitlepage
\let\maketitle\maketitlepage

\def\ifplaintex{\expandafter\ifx\csname documentclass\endcsname\relax}

\def\gtp{{\mathsurround=0pt\it $\cal G\mskip-2mu$eometry \&\ 
$\cal T\!\!$opology $\cal P\!$ublications}}  

\def\recd{{\small Received:\qua\receiveddate\ifx\reviseddate\relax
\else\qquad Revised:\qua\reviseddate\fi\par}} 

\def\lognumber#1{\def\thelognumber{#1}}
\def\volumenumber#1{\def\thevolumenumber{#1}}
\def\volumeyear#1{\def\thevolumeyear{#1}}
\def\papernumber#1{\def\thepapernumber{#1}}
\def\pagenumbers#1#2{\def\startpage{#1}\def\finishpage{#2}}
\def\published#1{\def\publishdate{#1}}

\def\received#1{\def\receiveddate{#1}}
\def\revised#1{\def\reviseddate{#1}}
\def\accepted#1{\def\accepteddate{#1}}

\long\def\asciiabstract#1{\long\def\theasciiabstract{#1}}

\let\\\par\let\thelognumber\relax\let\thevolumenumber\relax
\let\thepapernumber\relax\let\thevolumeyear\relax\let\startpage\relax
\let\finishpage\relax\let\publishdate\relax\let\receiveddate\relax
\let\reviseddate\relax\let\accepteddate\relax\let\theasciititle\relax
\let\theasciiauthors\relax
\let\theasciiabstract\relax

\let\theasciiemail\relax

\ifplaintex
\font\logobig=cmssbx10 scaled 3836
\font\logomed=cmssbx10 scaled 2557
\else
\font\logobig=cmssbx10 scaled 4200
\font\logomed=cmssbx10 scaled 2800
\fi

\long\def\makeagttitle{   
\count0=\startpage
\agt\hfill      
\hbox to 45truept{\vbox to 0pt{\vglue -13truept{\logomed A\kern -.37em{\logobig 
T}\kern -.38em G}\vss}\hss}
\break
{\small Volume \thevolumenumber\ (\thevolumeyear)
\startpage--\finishpage\nl
Published: \publishdate}

\vglue .25truein

{\parskip=0pt\leftskip 0pt plus
1fil\def\\{\par\smallskip}{\Large\bf\thetitle}\par\medskip} \vglue
0.05truein

{\parskip=0pt\leftskip 0pt plus 1fil\def\\{\par}{\sc\theauthors}
\par\medskip}%
 
\vglue 0.03truein 

{\small\leftskip 25truept\rightskip 25truept{\bf Abstract}\stdspace\theabstract

{\bf AMS Classification}\stdspace\theprimaryclass
\ifx\thesecondaryclass\relax\else; \thesecondaryclass\fi\par
{\bf Keywords}\stdspace \thekeywords\par}\vglue 7truept

}   

\ifplaintex
\font\phead=cmsl9 scaled 950
\font\pnum=cmbx10 scaled 913
\font\pfoot=cmsl9 scaled 950
\headline{\vbox to 0pt{\vskip -4.5mm\line{\small\phead\ifnum
\count0=\startpage ISSN 1472-2739 (on-line) 1472-2747 (printed)
\hfill {\pnum\folio}\else\ifodd\count0\def\\{ }%
\ifx\theshorttitle\relax\thetitle\else\theshorttitle\fi\hfill{\pnum\folio}
\else\def\\{ and }{\pnum\folio}\hfill\ifx\theshortauthors\relax\theauthors
\else\theshortauthors\fi\fi\fi}\vss}}
\footline{\vbox to 0pt{\vglue 0mm\line{\small\pfoot\ifnum\count0=\startpage
\copyright\ \gtp\hfill\else
\agt, Volume \thevolumenumber\ (\thevolumeyear)\hfill\fi}\vss}}
\else
\font=cmsl9 scaled 1050
\font\lnum=cmbx10 
\font=cmsl9 scaled 1050
\makeatletter
\def\@oddhead{{\small\lhead\ifnum\count0=\startpage ISSN 1472-2739 
(on-line) 1472-2747 (printed)\hfill {\lnum\number\count0}\else\ifodd\count0
\def\\{ }\ifx\theshorttitle\relax \thetitle \else\theshorttitle\fi\hfill
{\lnum\number\count0}\else\def\\{ and }{\lnum\number\count0}
\hfill\ifx\theshortauthors\relax 
\theauthors\else\theshortauthors\fi\fi\fi}}\def\@evenhead{\@oddhead}
\def\@oddfoot{\small\lfoot\ifnum\count0=\startpage\copyright\ \gtp\hfill\else
\agt, Volume \thevolumenumber\ (\thevolumeyear)\hfill\fi}
\def\@evenfoot{\@oddfoot}
\makeatother
\fi
\let\maketitlepage\makeagttitle
\let\makeshorttitle\maketitlepage
\let\maketitle\maketitlepage

\newwrite\gtoutfile
\long\gdef\makeheadfile{  
{\def\\{, }\def\s{ }
\immediate\openout\gtoutfile head.xxx
\immediate\write\gtoutfile{To: math@arxiv.org}
\immediate\write\gtoutfile{Subject: put OR rep NNNNN:ppppp}
\immediate\write\gtoutfile{--text follows this line--}
\immediate\write\gtoutfile{Proxy-for: \ifx\theasciiauthors\relax
\theauthors\else\theasciiauthors\fi\s<\ifx\theasciiemail\relax\theemail\else\theasciiemail\fi>}
\immediate\write\gtoutfile{\noexpand\\}
\immediate\write\gtoutfile{Authors: \ifx\theasciiauthors\relax
\theauthors\else\theasciiauthors\fi}
{\def\\{ }\immediate\write\gtoutfile{Title: \ifx\theasciititle\relax
\thetitle\else\theasciititle\fi}}
\immediate\write\gtoutfile{Subj-class: GT or SG, GR etc}
\immediate\write\gtoutfile{MSC-class: \theprimaryclass\ifx\thesecondaryclass\relax\else, \thesecondaryclass\fi}
\immediate\write\gtoutfile{Journal-ref: Algebr. Geom. Topol. \thevolumenumber\s
(\thevolumeyear) \startpage-\finishpage}
\immediate\write\gtoutfile{Comments: Published by Algebraic and
Geometric Topology at}
\immediate\write\gtoutfile{\s\s\s  http://www.maths.warwick.ac.uk/agt/AGTVol\thevolumenumber/agt-\thevolumenumber-\thepapernumber.abs.html}
\immediate\write\gtoutfile{\noexpand\\}
\immediate\write\gtoutfile{}
\ifx\theasciiabstract\relax
\immediate\write\gtoutfile{\theabstract}\else
\immediate\write\gtoutfile{\theasciiabstract}\fi
\immediate\write\gtoutfile{}
\immediate\write\gtoutfile{\noexpand\\}
\immediate\write\gtoutfile{}
\immediate\closeout\gtoutfile}}  

\def\maketitlepage{\makeagttitle\makeheadfile}
\let\makeshorttitle\maketitlepage
\let\maketitle\maketitlepage

\def\ifplaintex{\expandafter\ifx\csname documentclass\endcsname\relax}

\def\gtp{{\mathsurround=0pt\it $\cal G\mskip-2mu$eometry \&\ 
$\cal T\!\!$opology $\cal P\!$ublications}}  

\def\recd{{\small Received:\qua\receiveddate\ifx\reviseddate\relax
\else\qquad Revised:\qua\reviseddate\fi\par}} 

\def\lognumber#1{\def\thelognumber{#1}}
\def\volumenumber#1{\def\thevolumenumber{#1}}
\def\volumeyear#1{\def\thevolumeyear{#1}}
\def\papernumber#1{\def\thepapernumber{#1}}
\def\pagenumbers#1#2{\def\startpage{#1}\def\finishpage{#2}}
\def\published#1{\def\publishdate{#1}}

\def\received#1{\def\receiveddate{#1}}
\def\revised#1{\def\reviseddate{#1}}
\def\accepted#1{\def\accepteddate{#1}}

\long\def\asciiabstract#1{\long\def\theasciiabstract{#1}}

\let\\\par\let\thelognumber\relax\let\thevolumenumber\relax
\let\thepapernumber\relax\let\thevolumeyear\relax\let\startpage\relax
\let\finishpage\relax\let\publishdate\relax\let\receiveddate\relax
\let\reviseddate\relax\let\accepteddate\relax\let\theasciititle\relax
\let\theasciiauthors\relax
\let\theasciiabstract\relax

\let\theasciiemail\relax

\ifplaintex
\font\logobig=cmssbx10 scaled 3836
\font\logomed=cmssbx10 scaled 2557
\else
\font\logobig=cmssbx10 scaled 4200
\font\logomed=cmssbx10 scaled 2800
\fi

\long\def\makeagttitle{   
\count0=\startpage
\agt\hfill      
\hbox to 45truept{\vbox to 0pt{\vglue -13truept{\logomed A\kern -.37em{\logobig 
T}\kern -.38em G}\vss}\hss}
\break
{\small Volume \thevolumenumber\ (\thevolumeyear)
\startpage--\finishpage\nl
Published: \publishdate}

\vglue .25truein

{\parskip=0pt\leftskip 0pt plus
1fil\def\\{\par\smallskip}{\Large\bf\thetitle}\par\medskip} \vglue
0.05truein

{\parskip=0pt\leftskip 0pt plus 1fil\def\\{\par}{\sc\theauthors}
\par\medskip}%
 
\vglue 0.03truein 

{\small\leftskip 25truept\rightskip 25truept{\bf Abstract}\stdspace\theabstract

{\bf AMS Classification}\stdspace\theprimaryclass
\ifx\thesecondaryclass\relax\else; \thesecondaryclass\fi\par
{\bf Keywords}\stdspace \thekeywords\par}\vglue 7truept

}   

\ifplaintex
\font\phead=cmsl9 scaled 950
\font\pnum=cmbx10 scaled 913
\font\pfoot=cmsl9 scaled 950
\headline{\vbox to 0pt{\vskip -4.5mm\line{\small\phead\ifnum
\count0=\startpage ISSN 1472-2739 (on-line) 1472-2747 (printed)
\hfill {\pnum\folio}\else\ifodd\count0\def\\{ }%
\ifx\theshorttitle\relax\thetitle\else\theshorttitle\fi\hfill{\pnum\folio}
\else\def\\{ and }{\pnum\folio}\hfill\ifx\theshortauthors\relax\theauthors
\else\theshortauthors\fi\fi\fi}\vss}}
\footline{\vbox to 0pt{\vglue 0mm\line{\small\pfoot\ifnum\count0=\startpage
\copyright\ \gtp\hfill\else
\agt, Volume \thevolumenumber\ (\thevolumeyear)\hfill\fi}\vss}}
\else
\font=cmsl9 scaled 1050
\font\lnum=cmbx10 
\font=cmsl9 scaled 1050
\makeatletter
\def\@oddhead{{\small\lhead\ifnum\count0=\startpage ISSN 1472-2739 
(on-line) 1472-2747 (printed)\hfill {\lnum\number\count0}\else\ifodd\count0
\def\\{ }\ifx\theshorttitle\relax \thetitle \else\theshorttitle\fi\hfill
{\lnum\number\count0}\else\def\\{ and }{\lnum\number\count0}
\hfill\ifx\theshortauthors\relax 
\theauthors\else\theshortauthors\fi\fi\fi}}\def\@evenhead{\@oddhead}
\def\@oddfoot{\small\lfoot\ifnum\count0=\startpage\copyright\ \gtp\hfill\else
\agt, Volume \thevolumenumber\ (\thevolumeyear)\hfill\fi}
\def\@evenfoot{\@oddfoot}
\makeatother
\fi
\let\maketitlepage\makeagttitle
\let\makeshorttitle\maketitlepage
\let\maketitle\maketitlepage

\newwrite\gtoutfile
\long\gdef\makeheadfile{  
{\def\\{, }\def\s{ }
\immediate\openout\gtoutfile head.xxx
\immediate\write\gtoutfile{To: math@arxiv.org}
\immediate\write\gtoutfile{Subject: put OR rep NNNNN:ppppp}
\immediate\write\gtoutfile{--text follows this line--}
\immediate\write\gtoutfile{Proxy-for: \ifx\theasciiauthors\relax
\theauthors\else\theasciiauthors\fi\s<\ifx\theasciiemail\relax\theemail\else\theasciiemail\fi>}
\immediate\write\gtoutfile{\noexpand\\}
\immediate\write\gtoutfile{Authors: \ifx\theasciiauthors\relax
\theauthors\else\theasciiauthors\fi}
{\def\\{ }\immediate\write\gtoutfile{Title: \ifx\theasciititle\relax
\thetitle\else\theasciititle\fi}}
\immediate\write\gtoutfile{Subj-class: GT or SG, GR etc}
\immediate\write\gtoutfile{MSC-class: \theprimaryclass\ifx\thesecondaryclass\relax\else, \thesecondaryclass\fi}
\immediate\write\gtoutfile{Journal-ref: Algebr. Geom. Topol. \thevolumenumber\s
(\thevolumeyear) \startpage-\finishpage}
\immediate\write\gtoutfile{Comments: Published by Algebraic and
Geometric Topology at}
\immediate\write\gtoutfile{\s\s\s  http://www.maths.warwick.ac.uk/agt/AGTVol\thevolumenumber/agt-\thevolumenumber-\thepapernumber.abs.html}
\immediate\write\gtoutfile{\noexpand\\}
\immediate\write\gtoutfile{}
\ifx\theasciiabstract\relax
\immediate\write\gtoutfile{\theabstract}\else
\immediate\write\gtoutfile{\theasciiabstract}\fi
\immediate\write\gtoutfile{}
\immediate\write\gtoutfile{\noexpand\\}
\immediate\write\gtoutfile{}
\immediate\closeout\gtoutfile}}  

\def\maketitlepage{\makeagttitle\makeheadfile}
\let\makeshorttitle\maketitlepage
\let\maketitle\maketitlepage

\def\ifplaintex{\expandafter\ifx\csname documentclass\endcsname\relax}

\def\gtp{{\mathsurround=0pt\it $\cal G\mskip-2mu$eometry \&\ 
$\cal T\!\!$opology $\cal P\!$ublications}}  

\def\recd{{\small Received:\qua\receiveddate\ifx\reviseddate\relax
\else\qquad Revised:\qua\reviseddate\fi\par}} 

\def\lognumber#1{\def\thelognumber{#1}}
\def\volumenumber#1{\def\thevolumenumber{#1}}
\def\volumeyear#1{\def\thevolumeyear{#1}}
\def\papernumber#1{\def\thepapernumber{#1}}
\def\pagenumbers#1#2{\def\startpage{#1}\def\finishpage{#2}}
\def\published#1{\def\publishdate{#1}}

\def\received#1{\def\receiveddate{#1}}
\def\revised#1{\def\reviseddate{#1}}
\def\accepted#1{\def\accepteddate{#1}}

\long\def\asciiabstract#1{\long\def\theasciiabstract{#1}}

\let\\\par\let\thelognumber\relax\let\thevolumenumber\relax
\let\thepapernumber\relax\let\thevolumeyear\relax\let\startpage\relax
\let\finishpage\relax\let\publishdate\relax\let\receiveddate\relax
\let\reviseddate\relax\let\accepteddate\relax\let\theasciititle\relax
\let\theasciiauthors\relax
\let\theasciiabstract\relax

\let\theasciiemail\relax

\ifplaintex
\font\logobig=cmssbx10 scaled 3836
\font\logomed=cmssbx10 scaled 2557
\else
\font\logobig=cmssbx10 scaled 4200
\font\logomed=cmssbx10 scaled 2800
\fi

\long\def\makeagttitle{   
\count0=\startpage
\agt\hfill      
\hbox to 45truept{\vbox to 0pt{\vglue -13truept{\logomed A\kern -.37em{\logobig 
T}\kern -.38em G}\vss}\hss}
\break
{\small Volume \thevolumenumber\ (\thevolumeyear)
\startpage--\finishpage\nl
Published: \publishdate}

\vglue .25truein

{\parskip=0pt\leftskip 0pt plus
1fil\def\\{\par\smallskip}{\Large\bf\thetitle}\par\medskip} \vglue
0.05truein

{\parskip=0pt\leftskip 0pt plus 1fil\def\\{\par}{\sc\theauthors}
\par\medskip}%
 
\vglue 0.03truein 

{\small\leftskip 25truept\rightskip 25truept{\bf Abstract}\stdspace\theabstract

{\bf AMS Classification}\stdspace\theprimaryclass
\ifx\thesecondaryclass\relax\else; \thesecondaryclass\fi\par
{\bf Keywords}\stdspace \thekeywords\par}\vglue 7truept

}   

\ifplaintex
\font\phead=cmsl9 scaled 950
\font\pnum=cmbx10 scaled 913
\font\pfoot=cmsl9 scaled 950
\headline{\vbox to 0pt{\vskip -4.5mm\line{\small\phead\ifnum
\count0=\startpage ISSN 1472-2739 (on-line) 1472-2747 (printed)
\hfill {\pnum\folio}\else\ifodd\count0\def\\{ }%
\ifx\theshorttitle\relax\thetitle\else\theshorttitle\fi\hfill{\pnum\folio}
\else\def\\{ and }{\pnum\folio}\hfill\ifx\theshortauthors\relax\theauthors
\else\theshortauthors\fi\fi\fi}\vss}}
\footline{\vbox to 0pt{\vglue 0mm\line{\small\pfoot\ifnum\count0=\startpage
\copyright\ \gtp\hfill\else
\agt, Volume \thevolumenumber\ (\thevolumeyear)\hfill\fi}\vss}}
\else
\font=cmsl9 scaled 1050
\font\lnum=cmbx10 
\font=cmsl9 scaled 1050
\makeatletter
\def\@oddhead{{\small\lhead\ifnum\count0=\startpage ISSN 1472-2739 
(on-line) 1472-2747 (printed)\hfill {\lnum\number\count0}\else\ifodd\count0
\def\\{ }\ifx\theshorttitle\relax \thetitle \else\theshorttitle\fi\hfill
{\lnum\number\count0}\else\def\\{ and }{\lnum\number\count0}
\hfill\ifx\theshortauthors\relax 
\theauthors\else\theshortauthors\fi\fi\fi}}\def\@evenhead{\@oddhead}
\def\@oddfoot{\small\lfoot\ifnum\count0=\startpage\copyright\ \gtp\hfill\else
\agt, Volume \thevolumenumber\ (\thevolumeyear)\hfill\fi}
\def\@evenfoot{\@oddfoot}
\makeatother
\fi
\let\maketitlepage\makeagttitle
\let\makeshorttitle\maketitlepage
\let\maketitle\maketitlepage

\newwrite\gtoutfile
\long\gdef\makeheadfile{  
{\def\\{, }\def\s{ }
\immediate\openout\gtoutfile head.xxx
\immediate\write\gtoutfile{To: math@arxiv.org}
\immediate\write\gtoutfile{Subject: put OR rep NNNNN:ppppp}
\immediate\write\gtoutfile{--text follows this line--}
\immediate\write\gtoutfile{Proxy-for: \ifx\theasciiauthors\relax
\theauthors\else\theasciiauthors\fi\s<\ifx\theasciiemail\relax\theemail\else\theasciiemail\fi>}
\immediate\write\gtoutfile{\noexpand\\}
\immediate\write\gtoutfile{Authors: \ifx\theasciiauthors\relax
\theauthors\else\theasciiauthors\fi}
{\def\\{ }\immediate\write\gtoutfile{Title: \ifx\theasciititle\relax
\thetitle\else\theasciititle\fi}}
\immediate\write\gtoutfile{Subj-class: GT or SG, GR etc}
\immediate\write\gtoutfile{MSC-class: \theprimaryclass\ifx\thesecondaryclass\relax\else, \thesecondaryclass\fi}
\immediate\write\gtoutfile{Journal-ref: Algebr. Geom. Topol. \thevolumenumber\s
(\thevolumeyear) \startpage-\finishpage}
\immediate\write\gtoutfile{Comments: Published by Algebraic and
Geometric Topology at}
\immediate\write\gtoutfile{\s\s\s  http://www.maths.warwick.ac.uk/agt/AGTVol\thevolumenumber/agt-\thevolumenumber-\thepapernumber.abs.html}
\immediate\write\gtoutfile{\noexpand\\}
\immediate\write\gtoutfile{}
\ifx\theasciiabstract\relax
\immediate\write\gtoutfile{\theabstract}\else
\immediate\write\gtoutfile{\theasciiabstract}\fi
\immediate\write\gtoutfile{}
\immediate\write\gtoutfile{\noexpand\\}
\immediate\write\gtoutfile{}
\immediate\closeout\gtoutfile}}  

\def\maketitlepage{\makeagttitle\makeheadfile}
\let\makeshorttitle\maketitlepage
\let\maketitle\maketitlepage